\documentstyle{article}
\def \E {{\bf E}}

 \def \I {{\hbox{Im\,}}}

 \def \T {{\hbox { Tr }}}

 \def \i {{\hbox{i}}}

\def \l {\lambda}

\def \d {{\hbox {d}}}

\def \vep {\varepsilon}

\def \la {\langle}
\def \ra {\rangle}

\def \i {{\hbox {i}}}

\begin{document}

\author{A.KHORUNZHY{\footnote{Address after September, 2001:
D\'epartement de Math\'ematiques,
Universit\'e de Versailles Saint-Quentin, Versailles, FRANCE}}\\
 Institute for Low Temperature Physics,
Kharkov, \textsc{Ukraine} and \\
University Paris-7, \textsc{France}
\and
W.KIRSCH \\
Institute of Mathematics, Ruhr-University Bochum,\\
Bochum, \textsc{Germany}}
\title{ON ASYMPTOTIC EXPANSIONS AND SCALES OF SPECTRAL UNIVERSALITY \\
IN BAND RANDOM MATRIX ENSEMBLES }
%\date{ }
\maketitle

\begin{abstract}
We consider real random symmetric $N\times N$ matrices
$H$ of the band-type form with characteristic length $b$.
The
  matrix entries $H(x,y), x\le y$ are independent Gaussian random variables
  and have the variance proportional to
$ u({x-y\over b})$,
where $u(t)$ vanishes at infinity.
We study  the resolvent
$G (z) = (H-z)^{-1}$, $\I z\neq 0$  in the limit
$1\ll b\ll N$ and obtain explicit expression
$S(z_1,z_2)$ for the leading term
 of the first correlation function of the normalized trace
$\langle G(z)\rangle = N^{-1} {\hbox{ Tr }}G(z)$.

We examine $S(\l_1 +\i 0,\l_2-\i 0)$
on the local scale $\l_1 - \l_2 = {r\over N}$
and show that its asymptotic behavior
is
determined by the rate of decay of $u(t)$.
In particular, if $u(t)$ decays exponentially, then
$S(r) \sim - C\, b^2 N^{-1} r^{-3/2}$.
This expression is universal in the sense that
the particular form of $u$ determines the value of  $C>0$ only.
Our results
agree with those detected in both numerical and
theoretical
physics studies of spectra of band random matrices.

\end{abstract}

\section{Problem, motivation and results}

Random matrices play an important role in various fields of
mathematics and physics. The
eigenvalue distribution of large matrices
was initially considered by E.Wigner
to model the statistical properties of energy spectrum of heavy
nuclei (see e.g. the collection of early papers \cite{P}).
Further
investigations have led to numerous applications of random matrices of
infinite dimensions in such branches of theoretical physics as statistical
mechanics of disordered spin systems, solid state physics, quantum chaos
theory, quantum field theory and others (see monographs and reviews
\cite{BIZ,CPV,MGZ,H}). In mathematics, the spectral theory of
random matrices has revealed deep links
 with the orthogonal polynomials,
integrable systems, representation theory, combinatorics, non-commutative
probability theory
and other theories \cite{BI,DIZ,So,VDN}.

In present paper we deal with the family of real symmetric random matrices
that can be referred to as the band-type one. In the simplest
case the matrices have zeros outside of a band around the
principal diagonal. Inside of this band they are assumed to be jointly
independent random variables. The limiting transition considered is that the
band width $b$ increases at the same time as the dimension of the matrix $n$
does.

There is a large number of papers devoted to the use of
random matrices of this type in models of quantum chaotical
systems (see, e.g. \cite{SVZ} and references therein).
In these studies, one of the central topic is
related with the transition  between fully developed
chaos and complete integrability. The crucial observation made numerically
\cite{CMI} and then supported in the welth of theoretical physics papers
(see, for example \cite{FM,Syl})
is that
the ratio $b^2/n$
is the critical one for the corresponding transition
in spectral properties of band random matrices.

On the rigorous level,
the eigenvalue distribution of $H^{(n,b)}$ has been studied
\cite{BMP,CG,MPK}.
It is shown that the limiting eigenvalue distribution
exists, is non-random and depends on
the parameter $\alpha = \lim_{n\to\infty} b/n$.
However, the role of the ratio $b^2/n$ has not been revealed there.
Recently, a series of papers appeared where the band random matrices
are studied in the context of the non-commutative probability  theory \cite{Gui,Sh}.
These studies  also deal with the limit
$n,b\to\infty$
such that $\alpha>0$.

In present paper we are concentrated on the case of $\alpha=0$
represented by the limit
$$
1\ll b\ll n
$$
and study the first correlation function
of the resolvent of band random matrices.
We show that the ratio $\beta = \lim_{n\to\infty} b^2/n$ naturally arises
when one considers  the leading term of this correlation function on the local scale.
This can be regarded as the support of the conjecture
that the local properties of spectra of band random matrices depend
on the value of $\beta$.

Let us describe our results in more details.
We consider  the ensemble
$\{H^{(n,b)}\}$ of random $N\times N$ matrices, $N=2n+1$ whose entries
$H^{(n,b)}(x,y)$ have the variance
proportional to $u({x-y\over b})$, where $u(t)\ge 0$ vanishes at infinity.
We consider the resolvent
$G^{(n,b)}(z)=\left(H^{(n,b)}-z\right) ^{-1}$ and study the asymptotic expansion
of the
correlation function
$$
C_{n,b}(z_1,z_2) =
\E f_{n,b}(z_1)f_{n,b}(z_2) -
\E f_{n,b}(z_1)
\E f_{n,b}(z_2),
$$
where we denoted $f_{n,b}(z) = N^{-1} \T G^{(n,b)}(z)$.
Keeping $ z_i$ far from the real axis,
we consider
the leading term $S(z_1,z_2)$ of this expansion and
find explicit expression for it.
This term $S(r_1+\i 0,r_2-\i 0)$ regarded
on the local scale $r_1-r_2 = r/N$ exhibits different behavior
depending on the rate of decay of the profile  function $u(t)$.

Our main conclusion is that if
 $u(t)\sim \left| t\right| ^{-1-\nu }$ as $t\to\infty$, then
the value $\nu = 2$ separates two major cases.
If $\nu \in (1,2)$, then the limit of $S(r)$ depends on $\nu$.
If $\nu \in (2,+\infty)$, then
$$
{1\over Nb} S(r) =
- {\hbox{const}} \cdot {\sqrt N\over b } \cdot {1\over \vert r\vert^{3/2}}(1+o(1))\, .
$$
These results  are in agreement with those predicted  in theoretical physics
studies.
In particular, the last expression for $S$
coincides with the Altshuler-Shklovski asymptotics of the spectral
correlation function (see e.g. \cite{MFD}).

The paper is organized as follows.
In Section 2 we determine the family of ensembles and present several
already known results that will be needed. In Section 3 we formulate
our main propositions and describe the scheme of their proofs. To
illustrate this scheme, we present a short proof of the Wigner semicircle
law for Gaussian Orthogonal Ensemble of random matrices.
In Section 4 we study the correlation function
$C_{n,b}(z_1,z_2)$
and obtain the explicit expression $S(z_1,z_2)$ for its leading term.
In Section 5 we study
the self-averaging property of $G(z)$
and prove auxiliary facts used in Section 4.
 Expressions derived in Section 4 are analyzed in
Section 6, where the asymptotic behavior
of $S(z_1,z_2)$
is studied.
In Section 7 we give a summary of our
observations.

\section{Band Random Matrices and Wigner Law}

\subsection{The ensemble}

Let us consider the family ${\cal A}=\left\{ a(x,y),\;x\le
y,\,x,y\in {\bf Z}\right\} $   of jointly independent random
variables determined on the same probability space. We  assume that they
have joint Gaussian (normal) distribution with properties
$$
{\bf E}\,a(x,y)=0,\quad {\bf E}\,\left[ a(x,y)\right] ^2= v(1+\delta_{xy}),
\eqno(2.1)
$$
where we denote by $\delta$  the Kronecker symbol;
$$
\delta_{xy} =
\cases{0, & if $x\neq y,$\cr
              1, & if $x=y.$\cr}
$$
Here and below ${\bf E}$
denotes the mathematical expectation with respect
to the measure generated by the family ${\cal A}$.

Let $u(t), t\in {\bf R}$ be  a piece-wise continuous
function $ u(t)= u(-t)\ge 0$ satisfying conditions
$$
\sup_{t\in {\bf R}}\left|  u(t)\right| = \bar u <\infty
\eqno\,(2.2)
$$
and
$$
\int_{{\bf R}} u(t)dt=1.
\eqno\,(2.3)
$$
For simplicity, we assume $u(t)$ to be monotone for $t\ge 0$.

Given real parameter  $b>0 $, we introduce an infinite  matrix $U^{(b)}$
$$
U^{(b)}(x,y)= \frac 1b u\left( \frac{x-y}b\right) , \quad x,y \in {\bf Z}.
$$
and
determine the ensemble $\left\{ H^{(n,b)}\right\} $
as the family of real symmetric
matrices of the form
$$
H^{(n,b)}(x,y)=a(x,y)\sqrt{U^{(b)}(x,y)},\quad x\le y,
\;\;\vert x\vert,\vert y\vert \le n,
\eqno (2.4)
$$
where $b\le N,\;N=2n+1$ and the square root is assumed to be positive.

Let us note that the matrix (2.1) has the really band form when $U^{(b)}$ is
constructed with the help of a function $ u$ having a finite support,
say
$$
u(t) = \cases{ 1, & if $t\in ( -{1\over 2},{1\over 2})$,\cr
0, & otherwise.\cr}
$$
In this case the band width is less than or equal to $2b+1$.
If $b=N$, then matrices (2.4) coincide with those
belonging to Gaussian Orthogonal Ensemble (GOE) \cite{M}.
This random matrix ensemble
is determined as the family $\left\{ A_N\right\} $ of real symmetric
matrices
$$
A_N(x,y)=\frac 1{\sqrt{N}}a(x,y),\;\;x,y=1,\dots ,N,
\eqno\,(2.5)
$$
with $\left\{ a(x,y)\right\} $ belonging to ${\cal A}$ (2.1). GOE
together with its Hermitian and quaternion versions plays the fundamental
role in the spectral theory of random matrices.

Random symmetric
matrices (2.5) with independent arbitrary distributed random variables $a(x,y)$
 satisfying
(2.1) is referred to as the Wigner ensemble of random matrices.
This random matrix ensemble considered first by E. Wigner \cite{W}
is extensively studied in a series of papers
(see e.g. \cite{So} and references therein).
In particular, in paper
\cite{KP2} the resolvent technique is developed to study
the spectral properties of the Wigner ensemble.
Actually, we follow a version of this technique, but restrict
ourself with more simple case of Gaussian random variables.
More general case of arbitrary distributed random variables
would make the computations much more cumbersome.
Let us repeat that the main task of this paper
is to study  the role of the ratio between
$b $ and $N$ with respect to the spectral
properties of random matrices.

Finally, it should be noted that we restrict ourself with the ensemble of
real symmetric matrices for the sake of simplicity also. All results can be
obtained by using essentially the same technique for the Hermitian
analogue of $H^{(n,b)}$.

\subsection{Limiting eigenvalue distribution}

Eigenvalue distribution of matrices $H^{(n,b)}$  is described by
the normalized eigenvalue counting function
$$
\sigma (\lambda ;H^{(n,b)})=\#\{\lambda _j^{(n,b)}\le \lambda \}N^{-1},
\eqno
\,(2.6)
$$
where $\lambda _1^{(n,b)}\le \dots \le \lambda _N^{(n,b)}$ are eigenvalues
of $H^{(n,b)}$ . We denote by $f_{n,b}(z),\,z\in {\bf C}$ the Stieltjes
transform of the measure given by (2.6);
$$
f_{n,b}(z)=\int_{-\infty }^\infty \frac{d\sigma _{n,b}(\lambda )}{\lambda -z}
,\;\,\,{\I}\,z\neq 0.
\eqno\,(2.7)
$$
Given a Stieltjes transform $f(z)$, one can restore corresponding measure $
d\sigma (\lambda )$ with the help of the inversion formula
(see e.g. \cite{Don}).

The limiting behavior of (2.7) as $n,b\rightarrow \infty $ was studied
in a series of papers \cite{BMP,CG,KLH,MPK}.
It was proved in \cite{MPK}
that $f_{n,b}(z)$ converges as $n,b\to\infty$ in probability
to a nonrandom function that depends on the ratio $\alpha = \lim b/N$;
$$
p-\lim _{n,b\rightarrow \infty }\;\,f_{n,b}(z)=w_\alpha (z).
\eqno\,(2.8)
$$
In particular, if $\alpha=0$, then the function
$w_0(z) \equiv w(z)$ satisfies equation
$$
w(z)=\frac 1{-z-v w(z)}.
\eqno\,(2.9)
$$
The solution of this equation is unique
in the class of functions satisfying condition
$$
\I w(z) \I z\ge 0
$$
and can be represented in the form
$w(z)=\int (\lambda -z)^{-1}d\sigma _w(\lambda )$, where $\sigma _w(\lambda )$
is the famous semicircle (or Wigner) distribution  \cite{W}
with
the density
$$
\rho _w(\lambda )=\sigma _w^{\prime }(\lambda )=\frac 1{2\pi v}
\cases{\sqrt{4v-\lambda^2},& if $|\lambda|^2\le 4v$,\cr
              0,& if $|\lambda|^2\ge 4v $.\cr}
\eqno\,(2.10)
$$

This density was obtained first by E. Wigner \cite{W}
for eigenvalues of random matrices of the "full" form (2.5)
and can be also obtained as the limit (2.8) with $\alpha = 1$
$$
\sigma _w(\lambda )=\lim_{N\rightarrow \infty }\sigma (\lambda ;A_N)
\eqno
\,(2.11)
$$

Thus, one gets the same eigenvalue distribution in the opposite limiting
transitions of narrow $\alpha = 0$ and wide $\alpha = 1$ band
widths. It is known that in the intermediate regime $ 0<\alpha<1$ the
limiting distribution differs from the semicircle (2.11) \cite{MPK}.
In present paper we concentrate ourself on the most interesting case
$\alpha =0$.

In present paper we always consider the
case of $1\ll b\ll n$.
As we have noted, the paper is aimed to detect the role
of the  parameter $\beta = \lim_{N\to\infty} b^2/N$.
To avoid technical problems, we restrict ourself with the range
$$
b=n^\chi \quad 1/3< \chi< 1.
\eqno
\,(2.12)
$$
We are convinced that our results are valid on the whole range
$0<\chi <1$.

\section{Main Propositions and Scheme of the Proof}

The resolvent
$$
G^{(n,b)}(z)=\left( H^{(n,b)}-zI\right) ^{-1},\;{\hbox{Im}}\,z\neq 0
$$
is  widely exploited in the spectral theory of operators.
Its normalized trace $\left\langle G^{(n,b)}(z)\right\rangle $
coincides with the Stieltjes transform $f_{n,b}(z)$ (2.7);
$$
\left\langle G^{(n,b)}(z)\right\rangle =\frac 1N\;{\hbox{Tr}}\,\left(
H^{(n,b)}-zI\right) ^{-1}=\frac 1N\sum_{j=1}^N\frac 1{\lambda
_j^{(n,b)}-z}.
$$

The results of this section are related with the asymptotic behavior of
 $\langle G^{(n,b)}(z)\rangle $ in the limit (2.12),
with $z\in \Lambda_\eta  $,
$$
\Lambda_\eta  =\left\{ z\in C:\;\left| {\hbox{Im\,}}z\right|
\ge \eta \right\} \quad {\hbox{with}}\,\,\, \eta = 2\sqrt{v}+1.
\eqno
\,(3.1)
$$

\subsection{Main technical  results}

Our first statement concerns the pointwise convergence of the
diagonal entries $G^{(n,b)} (x,x;z), \, |x|\le n$ of the resolvent.
Let us determine the set
$$
 B_L\equiv B_L(n,b)=\left\{ x\in {\bf Z}:\left| x\right| \le
n-bL\right\}.
\eqno\,(3.2)
$$

\vskip 0.5cm

\noindent {\bf Theorem 3.1}

\noindent {\it Given $\vep>0$, there exists a natural  $L$
such that
$$
\sup_{x\in B_L} |G^{(n,b)} (x,x;z) - w(z)|\le \vep, \quad \forall z\in \Lambda_\eta ,
\eqno (3.3)
$$
for sufficiently large $b,n$.}

\vskip 0.5cm

The result of Theorem 3.1 is interesting by itself.
We shall  use it hardly in the proof of
the following statement concerning the correlation function
$$
C_{n,b}(z_1,z_2)={\bf E}\left\langle G(z_1)\right\rangle \left\langle
G(z_2)\right\rangle -{\bf E}\left\langle G(z_1)\right\rangle \,
{\bf E}\left\langle G(z_2)\right\rangle.
\eqno (3.4)
$$
\vskip 0.5cm

\noindent {\bf Theorem 3.2.}

\noindent {\it If $z_i\in \Lambda_\eta , i=1,2$, then in the limit
$n,b\to\infty$ (2.12)}

$$
C_{n,b}(z_1,z_2)=\frac 1{Nb}S(z_1,z_2)+o(\frac 1{Nb}).
\eqno \,(3.5)
$$
{\it The explicit term of $S(z_1,z_2)$ is given by relation}
$$
S(z_1,z_2)=\frac{2v}{\left( 1-vw_1^2\right)
\left( 1-vw_2^2\right) }Q(z_1,z_2),
\eqno \,(3.6)
$$
{\it where } $w_j\equiv w(z_j),\,j=1,2$ {\it and }$Q(z_1,z_2)$
{\it is given by the formula}
$$
Q(z_1,z_2)=
{1\over 2\pi} \int_{{\bf R}}\frac{ w_1^2w_2^2 \tilde u_F(p)}{\left[ 1-v
w_1 w_2 \tilde u_F(p)\right] ^2}\d p,
$$
{\it where we denote by $\tilde u_F(p)$ the Fourier transform of $u$}
$$
\tilde u_F(p)=\int_{{\bf R}} u(t)
e^{ipt}\d t.
$$

\vskip 0.5cm

It should be noted that in the case of GOE (2.6)
relation (3.5) is valid with $b$ replaced by $N$ and
expression
 (3.6)  takes the following  form (see e.g. \cite{FMP,KKP})
$$
S_{{\hbox{\begin{tiny}GOE\end{tiny}}}}(z_1,z_2)=\frac{2v }
{\left( 1-vw_1^2\right)
\left( 1-vw_2^2\right) }\frac{w_1^2w_2^2}
{\left[ 1-vw_1w_2\right] ^2}.
\eqno\,(3.7)
$$
Let us briefly explain why  (3.6) and (3.7)
lead to different asymptotic expressions on the local
scale determined as
$$
z_1^{(N)} = \lambda + {r_1\over N} + \i 0, \,\,\,
 z_2^{(N)} = \lambda + {r_2\over N} - \i 0
\eqno (3.8)
$$
with $\lambda \in {\hbox{ supp }} \d \sigma_w$ (2.10).
It follows from equality (2.9) that
$$
\frac{w_1^2w_2^2}
{\left[ 1-vw_1w_2\right] ^2} = \left( {w_1 - w_2\over z_1 - z_2}\right)^2
\eqno (3.9)
$$
This expression tends to infinity in the limit (3.8) and
$v w(z_1) w(z_2) \to 1$ as well.
But after dividing by $N^2$, one obtains from (3.7) and (3.9) that
$$
{1\over N^2 } S_{{\hbox{\begin{tiny}GOE\end{tiny}}}}(z_1, z_2) =
- {1 \over (r_1-r_2)^2} (1+o(1)).
\eqno (3.10)
$$
The left-hand side of (3.1) is usually called
the wide (or smoothed) version of the eigenvalue density correlation function
and the expression in the right-hand side of (3.10)
is derived by various methods for different
random matrix ensembles \cite{D,FMP,BZ,KKP}.

In Section 6 we study  $S(z_1,z_2)$  with the spectral parameters $z_1, z_2$
given by  (3.8).
Now the singularity
of $Q(z_1,z_2)$  is determined  by convergence of $1-v w_1w_2 \tilde u_F(p)$ to zero.
This convergence depends on the behavior of $\tilde u_F(p)$ around the origin $p=0$;
that is why the rate of decay of $u(t)$ at infinity
dictates the form of the limiting expression for $S$ in the local scale.

\subsection{The method and short proof of semicircle law}

We prove Theorem 3.1  in Sections 4 and 5. We are based on
the moment relations approach for resolvents of
random matrices proposed and developed in \cite{KP1,KP2,PF}.
This technique is proved to be rather general, powerful
and applicable to various random matrix ensembles.
We use a modified version of this approach needed to study rather
complex case of band random matrices.
To make the subsequent exposition more transparent,
let us describe the principal
points of this method in application to
the simplest case represented by GOE (2.6).

\subsubsection{Families of averaged moments}

In the early 70s F.Berezin observed \cite{Ber}
that regarding correlation functions
of the formal density of states $\;\rho _N(\theta )=\sigma _N^{\prime
}(\theta )$
$$
P_k^{(N)}(\Theta _k)={\bf E}\,\rho _N(\theta _1)\cdots \rho _N(\theta _k),
$$
$\Theta _k=(\theta _{1,}\dots ,\theta _k),$ one can derive for them a system
of relations that resembles equalities for correlation functions of
statistical mechanics. In this system $P_k^{(N)}$ is expressed via sum of $
P_{k-1}^{(N)},$ $P_{k+1}^{(N)}$ and some terms that vanish in the limit $
N\rightarrow \infty $. This can be rewritten in the vector form
$$
\vec P^{(N)}=\vec P_0+B\vec P^{(N)}+\vec \phi ^{(N)},
$$
with certain operator $B$ and vector $\vec \phi$ such that
 $\left\| B\right\| <1$ and $\left\| \Phi ^{(N)}\right\| =o(1)$ in
appropriate Banach space. These properties prove existence of $\lim
{}_{N\rightarrow \infty }\vec P^{(N)}=\vec P$;
the special form
of $B$ implies that the limiting $\vec P$ is nonrandom with the
components $\prod \rho(\theta_j)$.

This approach has got its rigorous formulation
on the base of the  resolvent approach used first in
the random matrix theory in the pioneering work
\cite{MP}.
Regarding the resolvent
$G_N$, the main subject is goven by the infinite system
of moments
$$
L_k^{(N)}(X_k,Y_k;Z_k)={\bf E}\,\prod_{j=1}^kG_N(x_j,y_j;z_j),\;k\in {\bf N}.
\eqno\,(3.11)
$$
The technique proposed in \cite{KP1,PF} and developed in
\cite{KP2} has been employed
in the study of eigenvalue distribution of various ensembles of random
operators and random matrices \cite{KKP,MPK}.

In present paper we use the moment relations approach in its modified
version.
The main observation here is that often it is
sufficient to study asymptotic behavior of $L_1^{(N)}$ and $L_2^{(N)}$
instead of the whole infinite family of the moments (3.11). This
considerably reduces amount of computations and makes the proofs more
transparent. To explain the principal steps of the proofs of Theorems 3.1
and 3.2, let us present here the short proof of the semicircle law for GOE.

\subsubsection{Derivation of system of relations}

The main ingredients in the derivation of moment relations are the resolvent
identity
$$
G^{\prime }(z)-G(z)=-G(z)\left( H^{\prime }-H\right) G^\prime (z),
\eqno\,(3.12)
$$
where $G^{\prime }(z)=(H^{\prime }-z)^{-1}$ , $G(z)=(H-z)^{-1}$ and equality
$$
{\bf E}\,\gamma F(\gamma )={\bf E}\,\gamma ^2\;{\bf E}\,F^{\prime }(\gamma ),
\eqno\,(3.13)
$$
where $\gamma $ is the Gaussian random variable with zero mathematical
expectation and $F(t),\,t\in {\bf R}$ is a nonrandom function such that all
integrals in (3.13) exist. Equality (3.13) is a simple consequence of the
integration by parts formula.

Let us consider (3.12) with $H^{\prime }=A_N$ (2.6) and $H=0$. We obtain
relation
$$
G_N(x,y)=\zeta \delta _{xy}-\zeta \sum_{s=1}^NG_N(x,s)A_N(s,y),\;\zeta
\equiv -z^{-1}.\eqno
\,(3.14)
$$
Regarding the normalized trace
$$
f_N(z)=\,N^{-1}\sum G_N(x,x)\equiv \left\langle
G_N\right\rangle
$$
and using (3.13), we obtain relation
$$
\E f_N=\zeta -\zeta \frac{v}{N^2}\sum_{x,s=1}^N(1+\delta _{xs}){\bf \,E}\,
\frac{\partial G_N(x,s)}{\partial A_N(s,x)}.
\eqno\,(3.15)
$$
One can easily find the partial derivatives with the help of (3.12).
Remembering that $H$ are real symmetric matrices, we have
$$
\frac{\partial G(x,y)}{\partial H(s,t)}=-\left[
G(x,s)G(t,y)+G(x,t)G(s,y)\right] \left( 1+\delta _{st}\right) ^{-1}.\eqno
\,(3.16)
$$
Substituting (3.16) into (3.15), we obtain the first main relation for
$L_1^{(N)}$
$$
\E f_N=\zeta +\zeta v [\E f_N]^2+\phi _1^{(N)}+\psi _1^{(N)},
\eqno\,(3.17)
$$
where
$$
\phi _1^{(N)}=\zeta v N^{-1}{\bf E}\,\left\langle G_N^2\right\rangle ,\,\,\,
{\hbox{and\ \,\,}}\psi _1^{(N)}=\zeta v {\bf E}\left\{
f_N^{\circ } f_N^{\circ }\right\}
$$
and we denoted by $\xi ^{\circ }$ the centered random variable
$$
\xi ^{\circ }=\xi -{\bf E}\,\xi .
$$
Clearly,
$\left\langle G\right\rangle ^{\circ }=\left\langle G^{\circ
}\right\rangle $ (here and till the end of the subsection we omit the
subscript $N$ in $G_N$).
If one can show that two last terms in (3.17) vanish as $N\to\infty$,
then convergence $\E f_N (z) \to w(z)$ will be proved.

We estimate the term $\phi_1$ with the help of two elementary
inequalities that hold for the resolvent of a real symmetric matrix:
$$
\left| f_N(z) \right| \le \left\| G(z)\right\| \le
 \left| {\hbox{Im}}\,z\right|^{-1}
$$
and
$$
\left\| G^2(z)\right\| \le \vert \I z \vert^{-2}.
$$
The last estimate implies that
$$
\sum_{s} \vert G(x,s)\vert^2 = \Vert G\vec e_x\Vert^2 \le \vert \I z\vert ^{-2}.
\eqno (3.18)
$$
Inequality (3.18) means that if $z\in \Lambda_\eta $, then $\left| \phi
_1^{(N)}\right| \le v\eta ^{-3}N^{-1}$.

\subsubsection{Selfaveraging property}

To show that $\lim_{N\to\infty} \psi_1^{(N)}=0$,
we prove that
the variance of $f_N$ vanishes
$$
{\bf Var} f_N =
{\bf E}\left| f_N^{\circ }\right| ^2=O(N^{-2}).
\eqno\,(3.19)
$$
It is clear that
$$
{\bf Var} f_N = {\bf E} \bar f_N^{\circ } f_N^{\circ}
={\bf E}\bar f_N^{\circ }f_N,
$$
where we denoted  $\bar f_N = f_N(\bar z)$.
Applying (3.14)
to the last factor $f_N$, we see that
$$
{\bf E}\bar f_N^{\circ }f_N = -{\zeta\over N}
\sum_{s,t=1}^N \E \left\{f_N^{\circ} G(x,s) A_N(s,x)\right\}.
$$
The using (3.13) and (3.16), we derive relation
$$
{\bf E}
 \bar f_N^{\circ } f_N =
 \zeta v\,{\bf E}\bar f_N^{\circ }f_N f_N +\phi_2^{(N)}+\psi _2^{(N)},
\eqno
\,(3.20)
$$
where
$
\phi_2^{(N)}=\zeta vN^{-1}{\bf E}\bar f_N^{\circ } \left\langle G^2
\right\rangle
$
and
$$
\psi _2^{(N)}=2\zeta vN^{-2}
{\bf E}\left\langle \bar G^2G\right\rangle .
$$
The useful observation is that
$$
{\bf E} \bar f_N^{\circ } f_N f_N
 ={\bf E} \bar f_N^{\circ}f_N {\bf E}f_N +
 {\bf E}\bar f_N^{\circ } f_N^{\circ }f_N .
 \eqno
\,(3.21)
$$
Using this identity and taking into account estimates (3.19), we derive from
(3.18) that
$$
{\bf E} \bar f_N^{\circ }f_N
\le v\eta ^{-1}\left\{ {\bf E} \bar f_N^{\circ }f_N\cdot {\bf E}
\left| f_N
\right| +{\bf E} \bar f_N^{\circ } f_N^{\circ } \cdot
\left| f_N \right|
\right\} +
$$
$$
v\eta ^{-3}N^{-1}{\bf E}\left| f_N^{\circ }
\right| +2v\eta ^{-4}N^{-2}.
$$
Taking into account that
$\left|  f_N^{\circ }
\right| ^2\equiv  \bar f_N^{\circ }
f_N^{\circ }$, we finally obtain
$$
{\bf E}\left| f_N^{\circ } \right| ^2\equiv {\bf E}
 \bar f_N^{\circ }f_N \le
$$
$$
 2v\eta ^{-2}{\bf E}\left| f_N^{\circ } \right|^2+
 v\eta ^{-3}N^{-1}\left( {\bf E}\left| f_N^{\circ}\right| ^2\right) ^{1/2}+
 2v\eta ^{-4}N^{-2}.
 \eqno\,(3.22)
$$
This immediately implies (3.19) provided $z\in \Lambda_\eta$ (3.1).
Obviously, $\psi_1^{(N)}$ admits the same estimate.

\subsubsection{The semicircle law and further corrections}

Returning back to (3.17) and gathering estimates for $\phi _1^{(N)}$ and $
\psi _1^{(N)}$ , one can easily derive that if $z\in \Lambda_\eta $ (3.10, then
$\lim_{N\rightarrow \infty
}g_N(z)=w(z)$ , with $w(z)$ given by (2.10).
Convergence of the Stieltjes transforms implies convergence of the
corresponding measures.
Thus the  semicircle law is proved.

It should be noted that relation (3.21) can be transformed into
$$
{\bf E} \bar f_N^{\circ } f_N f_N =
 2 {\bf E} \bar f_N^{\circ}f_N {\bf E}f_N +
 {\bf E}\bar f_N^{\circ } f_N^{\circ }f_N^{\circ} .
$$
Substituting this into  (3.18), we see that
$$
{\bf Var} f_N ={1\over N^2} {2\zeta v \over 1- 2\zeta v  \E f_N}
 \E \la \bar G^2 G\ra +
{1\over 1-2\zeta v \E f_N} \left(
\phi_2^{(N)} +
\E \left\{ \bar f_N^\circ
f_N^\circ f_N^\circ\right\}\right).
\eqno (3.23)
$$

Using the resolvent identity
$$
G(z_1)G(z_2) = - {G(z_1) - G(z_2)\over z_1 - z_2}, \quad
G(z_i) = \left( H-z_i \right)^{-1}
\eqno (3.24)
$$
and convergence of $\E f_N(z)$,
one can easily  find the limiting expression for $\E \la \bar G^2 G\ra$.
If one assumes that two last terms in (3.23) are values of the order
$o(N^{-2})$, then one arrives at (3.7) (see e.g. \cite{KKP}
for more details).

%%%%%%%%%%%%%%%%%%%%%%%%%%%%%%%%%%%%%%%%%%%%%%%%%%%%%%%%%%%%%%%%%%%%%%%
%%%%%%%%%%%%%%%%%%%%%%%%%%%%%%%%%%%%%%%%%%%%%%%%%%%%%%%%%%%%%%%%%%%%%%%
%%%%%%%%%%%%%%%%%%%%%%%%%%%%%%%%%%%%%%%%%%%%%%%%%%%%%%%%%%%%%%%%%%%%%%%
%%%%%%%%%%%%%%%%%%%%%%%%%%%%%%%%%%%%%%%%%%%%%%%%%%%%%%%%%%%%%%%%%%%%%%%
%%%%%%%%%%%%%%%%%%%%%%%%%%%%%%%%%%%%%%%%%%%%%%%%%%%%%%%%%%%%%%%%%%%%%%%
%%%%%%%%%%%%%%%%%%%%%%%%%%%%%%%%%%%%%%%%%%%%%%%%%%%%%%%%%%%%%%%%%%%%%%%
%%%%%%%%%%%%%%%%%%%%%%%%%%%%%%%%%%%%%%%%%%%%%%%%%%%%%%%%%%%%%%%%%%%%%%%
%%%%%%%%%%%%%%%%%%%%%%%%%%%%%%%%%%%%%%%%%%%%%%%%%%%%%%%%%%%%%%%%%%%%%%%
%%%%%%%%%%%%%%%%%%%%%%%%%%%%%%%%%%%%%%%%%%%%%%%%%%%%%%%%%%%%%%%%%%%%%%%
%%%%%%%%%%%%%%%%%%%%%%%%%%%%%%%%%%%%%%%%%%%%%%%%%%%%%%%%%%%%%%%%%%%%%%%
%%%%%%%%%%%%%%%%%%%%%%%%%%%%%%%%%%%%%%%%%%%%%%%%%%%%%%%%%%%%%%%%%%%%%%%
%%%%%%%%%%%%%%%%%%%%%%%%%%%%%%%%%%%%%%%%%%%%%%%%%%%%%%%%%%%%%%%%%%%%%%%
%%%%%%%%%%%%%%%%%%%%%%%%%%%%%%%%%%%%%%%%%%%%%%%%%%%%%%%%%%%%%%%%%%%%%%%
%%%%%%%%%%%%%%%%%%%%%%%%%%%%%%%%%%%%%%%%%%%%%%%%%%%%%%%%%%%%%%%%%%%%%%%
%%%%%%%%%%%%%%%%%%%%%%%%%%%%%%%%%%%%%%%%%%%%%%%%%%%%%%%%%%%%%%%%%%%%%%%

\section{Correlation Function of the Resolvent}

Our approach is to
apply systematically the scheme of subsection 3.2.2 to
get the leading term of the
correlation function  $C^{(n,b)}(z_1,z_2)$ (3.4).
This term is expressed
via the limit of the $\lim \E \la G^{(n,b)}(z)\ra = w(z)$
but we have to prove the pointwise version of this convergence
given by Theorem 3.1.
This and other auxiliary propositions are addressed in subsections 4.1 and 4.2.
In subsection 4.3 we give the proof of Theorem 3.2 on the base of these
statements.
In
what follows, we omit super- and subscripts $n$ and $b$ and do not indicate
the limits of
summation when no confusion can arise.

\subsection{Proof of Theorem 3.1}

Using relations (3.12)-(3.14) with obvious changes and repeating
computations of subsection 3.2.2, we obtain relation
$$
{\bf E}G(x,x)=\zeta +\zeta v{\bf E}G(x,x)U_G(x)+
\zeta v\sum_{\left| s\right|
\le n}{\bf E}\left[ G(x,s)\right] ^2U(s,x),
\eqno (4.1)
$$
where
$$
U_G(x)=\sum_{\left| s\right| \le n}G(s,s)U(s,x) =
{1\over b} \sum_{\left| s\right| \le n}G(s,s)u({s-x\over b}).
$$
Let us denote the average  ${\bf E}G(x,x)$ by $g(x)$ and rewrite
(4.1) in the following form
$$
g(x)=\zeta +\zeta v^2g(x)\,U_g(x)+\frac 1b\Phi(x)+\Psi(x),
\eqno
\,(4.2)
$$
where we denoted (cf. (3.17))
$$
\Phi(x)=\zeta v\sum_{\left| s\right| \le n}{\bf E}\left[
G(x,s)\right] ^2 u(\frac{s-x}b)
\eqno\,(4.3)
$$
and
$$
\Psi(x)=\zeta v{\bf E\,}G^{\circ }(x,x)\,U_G^{\circ }(x).
\eqno (4.4)
$$

Let us consider the solution  $\left\{ r(x),\,\left| x\right| \le n\right\} $
of  equation
$$
r(x)=\zeta +\zeta vr(x)U_r(x),\;\left| x\right| \le n.
\eqno\,(4.5)
$$
Given $z\in \Lambda_\eta  $ (3.1), one can prove that the system of equations
(4.5) is uniquely solvable in the set of $N$-dimensional vectors $\left\{
\vec r\right\} $ such that
$$
\left\| \vec r\right\| _1=\sup_{\left| x\right| \le n}\left| r(x)\right|
\le 2\eta ^{-1},\;\eta =\left| {\hbox{Im}}\,z\right|
\eqno\,(4.6)
$$
(see Lemma 4.1 at the end of this section). Certainly, $r(x)$ depends on
particular values of $z,\,n$ and $b$, so in fact we use denotation
$r(x)=r_{n,b}(x;z)$.

The following statements concern the differences
$$
D_{n,b}(x;z)=g_{n,b}(x;z)-r_{n,b}(x;z),\quad
\;d_{n,b}(x;z)=r_{n,b}(x;z)-w(z),
$$
where $w(z)$ is given as a solution of (2.9).

\vskip 0.5cm

{\bf Proposition 4.1.}

{\it Given $\vep>0$ ,   there
exists  a number $L=L(\vep )$ such that for all sufficiently large }$b${\it
\ and }$n${\it \ satisfying (2.13)
$$
\sup_{x\in B_L}\left| d_{n,b}(x;z)\right| \le \vep,\quad z\in \Lambda_\eta,
\eqno (4.7)
$$
with $B_L$ given by (3.2).
}

\vskip 0.5cm

{\bf Proposition 4.2.}

{\it If $z\in \Lambda_\eta  $ (3.1) and (2.13) holds, then
$$
\sup_{|x|\le n} |D_{n,b}(x;z)| = o(1), \quad n,b\to\infty.
\eqno (4.8)
$$
}
\vskip 0.5cm

Theorem 3.1
follows from (4.7) and (4.8).
Under the same conditions one can find $L'\ge L$ such that
$$
\sup_{x\in B_{L'}} |{\zeta \over 1-\zeta v U_g(x)} - w(z)|\le 2\vep.
\eqno (4.9)
$$
Relation (4.9) follows from (3.3) added by (4.6), a priori estimate
$$
\sup_{|x|\le n} | g (x) |\le {1\over |\I z|},
\eqno (4.10)
$$
and observation that $L'$ has to satisfy condition
$u(L-L')\le  \vep$.

\vskip 0.5cm

{\it Proof of Proposition 4.1}

Let us consider the constant function $w_x(z)\equiv w(z)$ satisfying (2.10)
that we rewrite in the following form similar to (4.5)
$$
w_x(z)=\zeta +\zeta vw_x(z)\frac 1b\sum_{\left| t\right| \le n}b
\delta
_{xt}w_t(z),\;\left| x\right| \le n.
$$
Subtracting this equality from (4.5), we derive that
$d(x)\equiv
d_{n,b}(x;z) $ verifies equality
$$
d(x)=\zeta vd(x)U_r(x)+\zeta vw(z)U_d(x)+\zeta vw^2(z)\left[ P_b+T(x)\right] ,
$$
where
$$
P_b=\frac 1b\sum_{t\in {\bf Z}} u\left( \frac tb\right) -\int_{-\infty
}^\infty u(s)ds\,
\eqno \,(4.11)
$$
and
$$
T_{n,b}(x)=\frac 1b\sum_{\left| t\right| \le n} u\left( \frac{t-x}
b\right) -\frac 1b\sum_{t\in {\bf Z}} u\left( \frac tb\right) .
\eqno (4.12)
$$

It is clear that  $\left| P_b\right| =o(1)$
as $b\rightarrow \infty $. Indeed, one can determine an even step-like function
$u_d(t), \, t\in {\bf R}$ such that
$$
u_d(t) = \sum_{k\in {\bf N} } u({k\over b})
{\hbox{I}}_{({k-1\over b}, {k\over b})}(t), \, t\ge 0.
$$
Then $u_d(t)\le u(t)$ and $u_d(t)\to u(t)$
as $b\to\infty$ and  the Beppo-L\'evy theorem
implies convergence of the corresponding integrals
of (4.10).

Taking into account equality
$$
r(x)={\zeta \over  1-v\zeta U_r(x)}\, ,
$$
we can  write that
$$
d(x)=vwr(x)U_d(x)+vw^2r(x)\left[ P_b+T_{n,b}(x)\right] ,\,
$$
where we denoted $w\equiv w(z)$. This relation, together with
estimates (4.6) and $|w(z)|\le |\I z|^{-1}$,
implies inequality
$$
\sup _{x\in B_L}\left| d(x)\right| \le \tau \left( \sup _{x\in
B_{L-1}}\left| d(x)\right| +\sup _{x\in B_L}\left| T_{n,b}(x)\right|+P_b\right) \le
$$
$$
\tau \sum_{j=0}^L\tau ^j\,\left( \sup _{x\in B_{L-j}}\left|
T_{n,b}(x)\right| +P_b\right) +\tau ^L\sup _{\left| x\right| \le n}\left|
d(x)\right| ,
\eqno \,(4.13)
$$
where $\tau \le v\eta ^{-2}<1.$ It is clear that due to monotonicity of
$u(t)$, one gets
$$
\sup _{x\in B_{L+1}}\left| T_{n,b}(x)\right| \le \sup _{x\in B_L}\left|
T_{n,b}(x)\right| \le \frac 2b\sum_{t=n-Lb}^\infty u\left( \frac
tb\right) \le 2\int_L^\infty u(s)ds.
$$
Given $\epsilon $ , one can find such a number $k$ that
$\tau ^k<\epsilon /4$. Than we derive from (4.13) that
$$
\sup _{x\in B_L}\left| d(x)\right| \le \tau \sum_{j=0}^k\tau ^j\,\sup _{x\in
B_{L-j}}\left| T_{n,b}(x)\right| +2\tau P_b+\epsilon /4\le
$$
$$
\tau \,\left( k+1\right) \sup _{x\in B_{L-k}}\left| T_{n,b}(x)\right| +2\tau
P_b+\epsilon /4.
$$
Now it is clear that (4.9) holds for sufficiently large $L$ and $b$.

Proposition 4.1 is proved.

\vskip 0.5cm

{\it Proof of Proposition 4.2.}

Subtracting (4.5) from (4.2), we obtain relation for $D(x)=D_{n,b}(x)$
$$
D(x)=\zeta vD(x)U_r(x)+\zeta vg(x)U_D(x)+\zeta v\left[ \frac 1b
\Phi(x)+\Psi(x)\right]
$$
that can be rewritten in the form
$$
D(x)=vg(x)r(x)U_D(x)+vr(x)\left[ \frac 1b\Phi(x)+\Psi(x)\right]
$$
Regarding this relation as the coordinate form of a vector equality,
one can write
that
$$
\vec D=v\left( I- W^{(g,r)}\right) ^{-1}\left[ \frac 1b\vec
\Phi^{(r)}+\vec \Psi^{(r)}\right] ,
$$
where we denote by  $W^{(g,r)}$ a linear operator acting on vectors $e$ with
components $e(x)$ as
$$
\left[ W^{(g,r)}e\right] (x)=vg(x)r(x)\sum_{|s|\le n} e(s) U(s,x)
$$
and vectors
$$
\vec \Phi _{n,b}^{(r)}(x)=r(x)\phi _{n,b}(x),\;\vec \Psi^{(r)}(x)=r(x)\Psi(x).
$$
It is easy to see that if $z\in \Lambda_\eta  $, then the  estimates
(4.6) and (4.10)
 imply inequality
$$
\left\| W^{(g,r)}\right\| \le {v\over \eta^2} < 1/2.
\eqno (4.14)
$$
Thus, to prove
Proposition 4.2, it is sufficient to show that
$$
\sup_x\left| \sum_s{\bf E}\left[ G(x,s)\right] ^2u({s-x\over b})\right|
=O(1),\;z\in \Lambda_\eta
\eqno (4.15)
$$
and
$$
\sup_x{\bf E}\left| G^{\circ }(x,x)\,U_G^{\circ }(x,x)\right| =
o\left(1\right) ,\;z\in \Lambda_\eta  .
\eqno\,(4.16)
$$
Relation (4.15) is a consequence of the bound (2.2) and inequality (3.18).
Relation
(4.16) reflects the selfaveraging property of $G^{(n,b)}$.
This question is addressed in the next subsection.
It should be noted that (4.16) will be proved independently from
computations of this subsection.
Assuming that this is done, we can say that Theorem 3.1 is proved.
We complete this subsection with the proof of the following auxiliary statement.

\vskip 0.5cm
{\bf Lemma 4.1.}

{\it Equation (4.5) has a unique solution in the class of vectors satisfying
condition (4.6).}

\vskip 0.3cm
{\it Proof.}

Let us consider the sequence of $N$-dimensional vectors $\left\{ \vec
r^{(k)},\;k\in {\bf N}\right\} $ determined by relations for their
components
$$
r^{(k+1)}(x)=\zeta +\zeta vr^{(k)}(x)U_{r^{(k)}}(x),\;r^{(1)}(x)=\zeta
,\;\left| x\right| \le n.
$$
Then it is easy to derive that if $\vec r^{(k)}$ satisfies (4.6) and $z\in
\Lambda_\eta  $ (3.1), then $\vec r^{(k+1)}$ also satisfies (4.6). The difference $
\chi _{k+1}(x)=r^{(k+1)}(x)-r^{(k)}(x)$ satisfies relations
$$
\chi _{k+1}(x)=\zeta v\chi _k(x)U_{r^{(k)}}(x)+\zeta vr^{(k-1)}(x)U_{\chi
_k}(x).
$$
Obviously, $\left\| \chi _{k+1}\right\| _1\le \alpha \left\| \chi _k\right\|
_1$ with $\alpha <1$ provided \thinspace $z\in \Lambda_\eta  $. Lemma is proved.$
\Box $

%%%%%%%%%%%%%%%%%%%%%%%%%%%%
%%%%%%%%%%%%%%%%%%%%%%%%%%

\subsection{The variance and selfaveraging property}

The asymptotic relation (4.15)
is a consequence of the  fact
that the variance of $G(x,x)$
$$
{\bf Var}\langle G^{(n,b)}\rangle ={\bf E}\left| \langle
G^{(n,b)}\rangle ^{\circ }\right| ^2
= \E \left\{ \la G^{(n,b)}(z)\ra^\circ \,\la G^{(n,b)}(\bar z)\ra^\circ\right\}
$$
vanishes as $n,b\to\infty$.
Instead of the direct proof of (4.15), we prefer to
present the whole  list of more general statements
needed in studies of the correlation function of $G$.
All of them can be proved independently of the Theorem 3.1
without use of its statement.

We start the list with
the following three relations that concern the moments of
diagonal elements of $G$.

\vskip 0.5cm
{\bf Proposition 4.3.}

{\it If $z\in \Lambda_\eta  $ (3.1), then the estimates}
$$
\sup_{\left| x\right| \le n}\;{\bf E}\left| G^{\circ }(x,x;z)\right|
^2=O(b^{-1}),
\eqno\,(4.17)
$$

$$
\sup_{\left| x\right| \le n}\;{\bf E}\left| U_G^{\circ }(x)\right|
^2=O(b^{-2}),
\eqno\,(4.18)
$$
{\it and}
$$
\sup_{\left| x\right| \le n}\;{\bf E}\left| U_G^{\circ }(x)\right|
^4=O(b^{-4}),
\eqno\,(4.19)
$$
{\it  hold. }

\vskip 0.5cm

The following statement concerns the mixed moments of variables $G^\circ(x,x;z)$
and their sums.
\vskip 0.5cm
{\bf Proposition 4.4.}

{\it If $z\in \Lambda_\eta$, then  relations}
$$
\sup_{\left| x\right| ,\left| y\right| \le n}\;\left| {\bf E}G^{\circ
}(x,x)U_G^{\circ }(y)\right| =O\left( b^{-2}\right) ,
\eqno\,(4.20)
$$
$$
\sup_{\left| x\right| \le n}\;\left| {\bf E}\left\langle G^{\circ
}\right\rangle G^{\circ }(x,x)\right| =O\left( n^{-1}b^{-1}+b^{-1}
\left[
{\bf Var}\left\langle G\right\rangle \right] ^{1/2}\right) ,
\eqno\,(4.21)
$$
{\it and}
$$
\sup_{\left| x\right| ,\left| y\right| \le n}\;\left| {\bf E}\left\langle
G^{\circ }\right\rangle G^{\circ }(x,x)U_G^{\circ }(y)\right| =O\left(
n^{-1}b^{-2}+b^{-2}\left[ {\bf Var}\left\langle G\right\rangle\right]
^{1/2}\right)
\eqno\,(4.22)
$$
{\it are true in the limit $n,b\to\infty$.}

\vskip 0.5cm

Finally, we formulate
\vskip 0.5cm
{\bf Proposition 4.5.}

{\it If $z\in \Lambda_\eta$, then relation}
$$
\sup_{\left| x\right| \le n}\left| {\bf E}\left\{ \left\langle G_1^{\circ
}\right\rangle \sum_s\left[ G_2(x,s)\right] ^2u_b^2(s,x)\right\} \right|
=O\left( n^{-1}b^{-2}+b^{-2}\left[ {\bf Var}\left\langle G
\right\rangle\right] ^{1/2}\right)
\eqno\,(4.23)
$$
{\it is true in the limit $n,b\to\infty$.}

\vskip 0.5cm
Let us not that the estimates (4.21)-(4.23) admit also the estimates
in terms of $n$ and $b$ that do not 
involve the variance of $\la G\ra$. 
However, derivation of the estimates would take more place and taime and we restrict
ourselves with the forms presented. 
It will be shown later that ${\bf Var}\left\langle G
\right\rangle = O(n^{-1}b^{-1})$. This fact together with the restriction
 (2.12) implies for (4.22) and (4.23) that 
$$
{1\over b^2} {1\over \sqrt {nb}} \ll {1\over nb}  
$$
that is sufficient for us. 
We prove Propositions 4.3-4.5 in Section 5.

\subsection{Toward the correlation function}

Let us have a more close look at the correlation function
$$
C_{n,b}(z_1,z_2) = \E \left\{\la G^{(n,b)}(z_1)\ra^\circ \,
\la G^{(n,b)}(z_2)\ra^\circ\right\}
$$
We follow the scheme described at the end of subsection 3.2 and
introduce variables $G_j(x,y)=G^{(n,b)}(x,y;z_j),\,j=1,2$.
To study the average
$$
{\bf E}\left\langle G_1^{\circ }\right\rangle
G_2(x,x)=R_{12}(x),
$$
we apply to $G_2(x,x)$ the resolvent identity (3.12)
and obtain relation
$$
R_{12}(x)=-\zeta _2\sum_{\left| s\right| \le n}{\bf E}\left\{ \left\langle
G_1^{\circ }\right\rangle G_2(x,s)a(s,x)\right\} \sqrt {U(s,x)},
$$
where $\zeta _2=-z_2^{-1}.$ We compute the last mathematical expectation
with the help of formulas (3.13) and (3.16) and obtain equality (cf. (4.1))
$$
R_{12}(x)=\zeta _2vR_{12}(x)U_{g_2}(x)+\zeta _2vU_{R_{12}}(x)g_2(x)+
$$
$$
2\zeta _2vN^{-1}\sum_s{\bf E}G_1^2(x,s)G_2(x,s)U(s,x)+\zeta _2v\left[
\Theta _{12}(x)+\Upsilon _{12}(x)\right],
$$
where we denoted $g_2(x)={\bf E}G(x,x;z_2)$,
$$
U_{g_{2}}(x)=\sum_{\left| s\right| \le n}g_{2}(s)U(s,x),
$$
$$
U_{R_{12}}(x)=\sum_{\left| s\right| \le n}R_{12}(s)U(s,x),
$$
$$
\Theta _{12}(x)={\bf E}\left\{ \left\langle G_1^{\circ }\right\rangle
\sum_{\left| s\right| \le n}\left[ G_2(x,s)\right] ^2 U(s,x)\right\} ,
$$
and
$$
\Upsilon _{12}(x)=E\left\{ \left\langle G_1^{\circ }\right\rangle
U_{G_2}^{\circ }(x)G_2^{\circ }(x)\right\} .
$$

Using denotation
$$
q_2(x) = {\zeta \over 1- \zeta v U_{g_2}(x)},
\eqno (4.24)
$$
we obtain the following relation for $R_{12}$
$$
R_{12}(x)= v q_2(x) g_2(x) U_{R_{12}}(x)+
{2v q_2(x)\over N} \sum_{s} F_{12}(x,s) U(s,x) +
$$
$$
v q_2(x) [\Theta _{12}(x)+\Upsilon _{12}(x)] ,
\eqno\,(4.25)
$$
where we denoted
$$
F_{12}(x,s) = \E G_1^2(x,s) G_2(x,s).
$$

The terms $\Theta$ and $\Upsilon$
can be estimated with the help of Propositions 4.3-4.5.
As we shall see in the next subsection, they do not contribute
to the leading term of $R_{12}$.
To obtain the explicit expression for the this leading term,
it is necessary to study  in detail the variable $F_{12}$.
Now let us formulate corresponding statement and
the auxiliary
relations
needed.

\vskip 0.5cm
\newpage

{\bf Proposition 4.6.}

{\it If $z\in \Lambda_\eta$, (3.1), then for arbitrary positive $\vep $
and large enough values of }$b${\it \ and }$n${\it \ (2.13) there exists
the set $B_L$ (3.2) with $L$ such that }
$$
\sup_{x\in B_L}\left| b[F_{12}U](x,x)-
{1\over 2\pi}
\frac{w_1^2w_2}
{1-v w_1^2}
\int_{{\bf R}}\frac{\tilde u_F(p)}{\left[ 1-vw_1 w_2
\tilde u_F(p)\right] ^2}\d p \right| \le \vep.
\eqno\,(4.26)
$$
\vskip 0.5cm

The proof of Proposition 4.6 is based on the similar statement
formulated for the product $G_1 G_2$.

\vskip 0.5cm

{\bf Proposition 4.7.}

{\it Given positive $\vep $, there exists such }$L$
{\it \ that relations }
$$
\sup_{x \in B_L}\left| b\sum_{\left| s\right| \le n}{\bf E}
G_1(x,s)G_2(x,s)\,U^k\left( s,x\right) -
{1\over 2\pi}
\int_{{\bf R}}
\frac{w_1w_2\tilde
u_F^k(p)}{1-vw_1w_2\tilde u_F(p)} \d p\right| \le \vep
\eqno\,(4.27)
$$
{\it and }
$$
\sup_{x\in B_L}\left| \sum_{\left| s\right| \le n}{\bf E}
G_1(x,s)G_2(x,s)\,-\frac{w_1w_2}{1-vw_1w_2}\right| \le \epsilon
\eqno (4.28)
$$
{\it hold for all $k\ge 1$, all $z_i\in \Lambda_\eta $ and
 large enough values of }$b${\it .}

\vskip 0.5cm

{\it Remark.}
In the case when $z_1\neq z_2$,  relation (4.28) can be
derived from
 the  resolvent identity
(3.24) with the help of the convergence (3.3) and the explicit form of
$w(z)$ (2.9).

\vskip 0.5cm

We prove Proposition 4.6 in the next subsection.
Relations (4.27) and (4.28) will be proved in Section 5.

\subsection{Proof of Proposition 4.6 and Theorem 3.2}

Let us assume that relations (4.27) and (4.28) are true and show that under
conditions of Theorem 3.2 the leading term of $R_{12}$ is of the order $
O(n^{-1}b^{-1})$ and terms $\Theta _{12}$ and $\Upsilon _{12}$ of (5.2) do
not contribute to it. We rewrite (4.25) in the form
$$
R_{12}(x)=vg_2(x)q_2(x)U_{R_{12}}(x)+2vq_2(x)N^{-1}\left[ F_{12}U\right]
(x,x)+
$$
$$
vq_2(x)\left[ \Theta _{12}(x)+\Upsilon _{12}(x)\right] .
\eqno\,(4.29)
$$
Let us denote $r_{12} = \sup_{|x|\le n} |R_{12}(x)|$.

Taking into account $U(x,y)\le \bar u/b$ (2.2) and
using inequalities of the form (3.19),
it is easy to see that if $z_i\in\Lambda_\eta$, then
$$
{1\over N}\left[ F_{12}U\right] (x)\le {1\over Nb}\E
\left( \sum_s |G_1^2(x,s)|^2\right)^{1/2} \left(\sum_s  |G_2(s,x)|^2\right)^{1/2}
= O\left( {1\over nb}\right) .
$$
Regarding this estimate and relations (4.22), (4.23)
we easily derive from (4.29) inequality (cf. (3.22))
$$
r_{12}\le
{v\over \eta^2 } r_{12} +
{C\over bn} +{1\over b^2} \sqrt{r_{12}}
$$
with some constant $C$.
Since $r_{12}$ is bounded for all $z\in \Lambda_\eta  $,
then
$$
r_{12}=O\left( {1\over nb}+{1\over b^4}\right) .
$$
Now condition (2.12) implies that
$
r_{12}=O(1/nb)$ and therefore
the general form of (3.5) is demonstrated.

Substituting  (3.5) into the estimates (4.22) and (4.23),
we obtain that
$$
\left\| \Theta _{12}\right\| _1=o\left( {1\over nb}\right) \,\,
{\hbox{ and}\,\;}
\left\| \Upsilon _{12}\right\| _1=o\left( {1\over nb}\right) .
$$
Thus, these terms of (4.29) do not contribute to the leading term
of $R_{12}$.
To find this term  in explicit form,
we need the result of Proposition 4.6.

\vskip 0.5cm

{\it Proof of Proposition 4.6.}

Regarding $F_{12}(x,y)={\bf E}
\,G_1^2(x,y)G_2(x,y)$, we
 apply to $G_2$ the resolvent identity (3.12).
Computations similar to those of subsection 3.2.2 lead us to equality
$$
F_{12}(x,y)=\zeta _2\delta _{xy}{\bf E}G_1^2(x,x)+\zeta _2v\left[
t_{12}U\right] (x,y)\;{\bf E}G_1^2(y,y)+
$$
$$
\zeta _2v\left\{ \left[ F_{12}U\right]
(x,y)\;g_1(y)+F_{12}(x,y)U_{[g_2]}(y)+\Gamma(x,y)\right\} ,
\eqno\,(4.30)
$$
where
$$
t_{12}(x,y)={\bf E}T_{12}(x,y)={\bf E}G_1(x,y)G_2(x,y),
$$
and the vanishing terms are denoted by $\Gamma = \Gamma_1 + \Gamma_2 + \Gamma_3$:
$$
\Gamma _1(x,y)={\bf \sum_sE}\left\{
G_1(x,y)G_1^2(s,y)G_2(x,s)+2G_1^2(x,y)G_2(s,y)G_2(x,s)\right\} U(s,y),
$$
$$
\Gamma _2(x,y)={\bf E}\left\{ \left[ T_{12}U\right] (x,y)\left[
G_1^2(y,y)\right] ^{\circ }\right\} +{\bf E}\left\{ \left[ F_{12}U\right]
(x,y)G_1^{\circ }(y,y)\right\} ,
$$
and
$$
\Gamma _3(x,y)={\bf E}F_{12}(x,y)U_{[G_2]}^{\circ }(y).
$$
Indeed, it is easy to show that
$$
\sup_{x,y} |\Gamma_j(x,y)| = O(b^{-1}), \quad z_1, z_2\in \Lambda_\eta .
\eqno (4.31)
$$
This can be done with the help of inequality (3.19) and
relations (4.17), (4.18), and (4.23).

Using definition of $q_2(x)$ (4.24),
we rewrite (4.30)  as
$$
F_{12}(x,y)=
vg_1(x)q_2(y)\left[ F_{12}U\right] (x,y)+
$$
$$
R^{(1)} (x,y) + R^{(2)}(x,y) +
v\tilde \Gamma(x,y),
\eqno\,(4.32)
$$
where we denoted
$$
R^{(1)}(x,y) = q_2(x)\,{\bf E}G_1^2(x,x)\delta _{xy},
\eqno (4.33a)
$$
$$
 R^{(2)}(x,y) = vq_2(y)\left[
t_{12}U\right] (x,y)\;{\bf E}G_1^2(y,y)
\eqno (4.33b)
$$
and $\tilde \Gamma(x,y) = \Gamma(x,y) q_2(y)$.
Let us note that $\vert R{(1)}\vert\le \eta^{-3}$ and
$\vert R^{(2)}\vert \le v \eta^{-5}$ for $z_i\in \Lambda_\eta$.

Let us determine the
linear operator $W$ that acts on $N\times N$ matrices $F$
according to the formula
$$
[WF](x,y)= v g_1(x)
\left[\sum_{\left| s\right| \le
n}F(x,s)U(s,y) \right]q_2(y).
$$
The a priori estimates
$
\vert g_1(x)\vert \le \vert \I z_1\vert^{-1},
$
and
$
\vert q_2(x)\vert \le \vert \I z_2\vert^{-1}
$
imply inequality (cf. (4.14))
$$
\Vert W\Vert_{(1,1)} \le  {v\over \eta^2}  <{1\over 2},\quad z_i\in \Lambda_\eta ,
\eqno (4.34)
$$
where the norm of $N\times N$ matrix $A$ is determined as
$\Vert A\Vert_{(1,1)} = \sup_{x,y} \vert A(x,y)\vert$.
This estimate is verified by direct computation of $\Vert W A\Vert_{(1,1)}$
with $\Vert A\Vert_{(1,1)} = 1$.

Then (4.32) can be rewritten as
$$
F_{12}(x,y)=v\sum_{m=0}^\infty
\left[W ^m
\left( R^{(1)} + R^{(2)}  +v\tilde \Gamma \right)\right] (x,y).
\eqno\,(4.35)
$$

The next steps  of the proof of  (4.26) are very elementary.
We consider the first $M$ terms of the infinite series
and use the decay of the matrix elements $U(x,y) = U^{(b)}(x,y)$ .
Indeed, if one considers (4.33) with
$x$ and $y$ taken far enough from the endpoints
$-n$, $n$, then the variables $g_1(s)$, $q_2(t)$ enters
into the finite series with $s$ and $t$ also far from the endpoints.
Then one can use relations (3.3) and (4.9)  and replace
$g_1$ and $q_2$ by the constant values $w_1$ and $w_2$, respectively.
This substitution leads simplifies expressions with
the error terms that vanish as $n,b\to\infty$.
The second step is similar.
It is to show that we can use Proposition  4.7
and replace the terms $R^{(1)}$ and $R^{(2)}$ of
the finite series of (4.33) by corresponding
expressions given by formulas (4.27) and (4.28).

Let us start to perform this program.
Taking   into account the estimate of $\Gamma$
and using  boundedness of terms $R^{(1)}$
and $R^{(2)}$, we can deduce from (4.35) that
$$
b\sum_s F_{12} (x,s) U(s,x) =
bv \sum_{m=0}^M \left[W^m (R^{(1)} + R^{(2)}) U\right](x,x) + \Delta^{(1)}(x,x),
\eqno (4.36)
$$
where $M$ is such that given $\vep>0$, $\vert \Delta^{(1)}(x,x)\vert <\vep$
for large enough $b$ and $n$.

Now let us find such $h$ that
the following holds
$$
u(t) \le \vep, \forall |t|\ge h, \quad {\hbox {and  }}
\int_{\vert t\vert \ge h} u(t) \d t \le \vep.
$$
We determine  matrix
$$
\hat U(x,y) =
\cases{ U(x,y), & if $ \vert x-y\vert \le bh$;\cr
0, & if $ \vert x-y\vert > bh$\cr}
$$
and denote by $\hat W$ corresponding linear operator
$$
[\hat W F](x,y) = v g_1(x)\left[  \sum_{\vert s \vert \le n} F(x,s)
 \hat U(s,y) \right]
q_2(y).
$$
Certainly, $\hat W$ admits the same estimate as (4.34).

Given $\vep>0$, let $L$ the largest number among those
required by conditions of  Propositions 4.1 and 4.7.
Let us denote by $Q$ the first natural greater than $(M+k)h$.
Then one can write that
$$
bv \sum_{m=0}^M \left[W^m (R^{(1)} + R^{(2)}) U\right](x,x)  =
$$
$$
bv \sum_{m=0}^M \left[(v\hat W)^m (R^{(1)} + R^{(2)}) \hat U\right](x,x)
+\Delta^{(2)}(x,x),
$$
where
$$
\sup_{x\in B_{L+Q}} \vert \Delta^{(2)}(x,x)\vert \le \vep, \quad {\hbox{as}}\,\,
n,b\to\infty.
\eqno (4.37)
$$
The proof of (4.37) uses  elementary computations.
Indeed, $ \Delta^{(2)}(x,x)$ is represented as the sum
of $M+1$  terms of the form
$$
bv^{m+1}\sum^*_{\vert x_i\vert \le n}
[g_1(x)]^m F(x,x_1) U(x,x_1) q_2(x_1) \cdots U(x_{m-2},x_{m-1})q_2(x_{m-1})
\times
$$
$$
[R^{(1)} + R^{(2)}](x_{m-1},x_{m})
U(x_{m}, x),
$$
where the sum is taken over the values of $x_i$ such that
 $\vert x_j - x_{j+1} \vert >bh$
at least for one of the numbers $j\le m$.

Now
remembering
the a priori bounds for $R^{(1)}$ and $R^{(2)}$, estimates like (4.13)
and
taking into account the diagonal form of $R^{(1)}$, one obtains
the following estimate of $\Delta^{(2)}$ by two terms
$$
\sup_{|x|\le n} \vert \Delta^{(2)}(x,x)\vert\le
\sum_{m=0}^M  {bv^{m+1} \over \eta^{2m+3}}\sum^*_{\vert x_i\vert \le n}
U(x,x_1) U(x_{1},x_{2})
\cdots U(x_{m}, x) +
$$
$$
\sum_{m=0}^M  {v^{m+1}\over \eta^{2m+5}}
\sum^*_{\vert x_i\vert \le n}
 U(x,x_1)  \cdots U(x_{m-2},x_{m-1})
U(x_{m}, x).
\eqno (4.38)
$$
Regarding the first term in the right-hand side of (4.38)
and assuming that $\vert x_j - x_{j+1}\vert > bh$,  one can observe that
for large enough $b,n$
$$
\sum_{\vert x_j\vert \le n} U^j(x,x_j) \vep U^{m-j} (x_{j+1},x) \le \vep.
$$
Indeed,
$$
\sum_{\vert x_i\vert \le n}
U(x,x_1) U(x_{1},x_{2})
\cdots U(x_{j-1}, x_j) \le
$$
$$
\sum_{ x_i\in {\bf Z}}
U(x,x_1) U(x_{1},x_{2})
\cdots U(x_{j-1}, x_j) \le
\left[ \int_{-\infty}^\infty u(t) \d t +{u(0)\over b}\right]^j\le
(1+\bar u/b)^j.
$$

Let us also mention here  that given $\vep>0$,
one has for large enough $n,b$ that
$$
\sup_{x\in B_{L+Q}} \left\vert \sum_{\vert s\vert \le n} U^j(x,s) -1
\right\vert\le \vep,
\eqno (4.39)
$$
where $j\le M$.
This follows from elementary computations
related with the differences (4.11) and (4.12) that vanish
in the limit $1\ll b\ll n$.

Similar but a little more modified reasoning
can be used to estimate the second
term in the right-hand side of (4.38).
Now one can write  that
$$
\sup_{|x|\le n} \vert \Delta^{(2)}(x,x)\vert\le 2\vep \sum_{m=0}^M
{m v^{m+1} \over \eta^{2m+2}} \le \vep.
$$

Regarding the right-hand side of (4.37)
with $x\in B_{L+Q} $, one observes that
the summations run over such values of $x_i$ that
$\vert x-x_1\vert \le bh,  \, \vert x_{i}-x_{i+1}\vert \le bh$,
and thus  $x_j\in B_L$ for all $j\le k+m-1$.
This means that we can apply relations (3.3) and (4.9)  to the right-hand
side of (4.37)
and
 replace $g_1(x) $ by $w(z_1)$,
$q_2(x)$ by $w(z_2)$. We derive from (4.36)
 that
$$
(F_{12}U^k)(x,x) = bv w_2\sum_{m=0}^M (vw_1 w_2)^m  (\hat U^m)(x,s)
\left[R^{(1)} + R^{(2)}\right](s,t) \hat U(t,x) + \Delta^{(3)}(x,x)
$$
with
$$
\sup_{x\in B_{L+Q}} \vert \Delta^{(3)}(x,x)\vert\le 4\vep.
$$
Finally, applying Proposition (4.7) to the expressions involved in $R$
and taking into account that
$$
\sup_{x\in B_{L+Q}} \vert b U^{m+1}(x,x) - {1\over 2\pi}\int \tilde u_F ^{m+1}(p) \d p
\vert \le \vep,
\eqno (4.40)
$$
 we obtain equality
$$
(F_{12}U)(x,x) = {v\over 2\pi} {w_1^2 w_2\over 1-vw_1^2} \sum_{m=0}^M (vw_1w_2)^m
\int \tilde u_F ^{m+1}(p) \d p +
$$
$$
{v\over 2\pi}{w_1^2 w_2\over 1-vw_1^2} \sum_{m=0}^M (vw_1w_2)^m
\int {\tilde u_F ^{m+1}(p) \over 1- v w_1 w_2 \tilde u_F} \d p +
\Delta ^{( 5)} (x,x)
\eqno (4.41)
$$
with
$$
\sup_{x\in B_{L+Q}}
\vert \Delta ^{( 5)} (x,x)\vert \le \vep \quad b,n\to\infty.
$$

Passing back in (4.41) to the infinite series and
simplifying them, we arrive
at the expression standing in the right-had side of (4.26).
Proposition is proved.
$\Box$

\vskip 0.5cm

Let us complete the proof of  Theorem 3.2.
Remembering estimate (4.14), we can
iterate relation (4.29) and obtain that
$$
R_{12}(x) = {2v q_2(x) \over Nb} \sum_{m=0}^\infty
[ (W^{(g_2,q_2)})^m \vec f_{12} ](x) + o(1/nb),
$$
where we denoted $\vec f_{12}(x) = b q_2(x) [F_{12} U](x,x)$.
Regarding the trace
$$
{1\over N} \sum_{\vert x\vert \le n} R_{12}(x) =
{1\over N} \sum_{x \in B_L} R_{12}(x)(1 + o(1)),
$$
and repeating the arguments of the proof of Proposition 4.6 presented
above,
we can
write
that
$$
R_{12}(x) = {2v w_2\over Nb} \sum_{m=0}^M
\sum_{t} (b F_{12} U)(t,t) (v  w_2^2  U)^m(t,x) + \Delta^{(6)}(x,x)
$$
with $\sup_{x\in B_L} \vert \Delta^{(6)}(x,x)\vert \le vep'$ provided
$n,b\to\infty$ (2.12).
Finally, observing that $(b F_{12} U)(t,t)$ asymptotically does not depend on $t$
(4.26), 
we arrive, with the help of (4.39),  at the expression (3.6). Theorem 3.2 is proved.

%%%%%%%%%%%%%%%%%%%%%%%%%%%%%%%%%%%%%%%%%%%%%%%%%%%%%%

%%%%%%%%%%%%%%%%%%%%%%%%%%%%%%%%%%%%%%%%%%%%%%%%%%%%%%
%%%%%%%%%%%%%%%%%%%%%%%%%%%%%%%%%%%%%%%%%%%%%%%%%%%%%%
%%%%%%%%%%%%%%%%%%%%%%%%%%%%%%%%%%%%%%%%%%%%%%%%%%%%%%
%%%%%%%%%%%%%%%%%%%%%%%%%%%%%%%%%%%%%%%%%%%%%%%%%%%%%%
%%%%%%%%%%%%%%%%%%%%%%%%%%%%%%%%%%%%%%%%%%%%%%%%%%%%%%
%%%%%%%%%%%%%%%%%%%%%%%%%%%%%%%%%%%%%%%%%%%%%%%%%%%%%%
%%%%%%%%%%%%%%%%%%%%%%%%%%%%%%%%%%%%%%%%%%%%%%%%%%%%%%
%%%%%%%%%%%%%%%%%%%%%%%%%%%%%%%%%%%%%%%%%%%%%%%%%%%%%%
%%%%%%%%%%%%%%%%%%%%%%%%%%%%%%%%%%%%%%%%%%%%%%%%%%%%%%
%%%%%%%%%%%%%%%%%%%%%%%%%%%%%%%%%%%%%%%%%%%%%%%%%%%%%%
%%%%%%%%%%%%%%%%%%%%%%%%%%%%%%%%%%%%%%%%%%%%%%%%%%%%%%

\section{Proof of auxiliary statements}

{\it Proof of Proposition 4.3}

Let us consider the average ${\bf E}\,G_1^{\circ }(x,x)G_2(y,y)$ and derive
for it, with the help of formulas (3.12), (3.13) and (3.16)
equality
$$
{\bf E}\,G_1^{\circ }(x,x)G_2(y,y)=\zeta _2v{\bf E}\,G_1^{\circ
}(x,x)G_2(y,y)U_{G_2}(y)+
$$
$$
\zeta _2v\sum_s{\bf E}\,G_1^{\circ }(x,x)\left[ G_2(y,s)\right]^2
U(s,y)+
$$
$$
2\zeta _2v\sum_s{\bf E}\,G_1(x,s)G_1(y,x)G_2(y,s) U(s,y).
$$

Applying to the first term of this equality
the analogue of identity (3.21) and using  $q_2(x)$ (4.24), we obtain that
$$
{\bf E}\,G_1^{\circ }(x,x)G_2(y,y)=v q_2(y){\bf E}\,G_1^{\circ
}(x,x)G_2(y,y)U^\circ_{G_2}(y)+
$$
$$
vq_2(y) \sum_s{\bf E}\,G_1^{\circ }(x,x)\left[ G_2(y,s)\right]^2
U(s,y)+
$$
$$
2v q_2(y) \sum_s{\bf E}\,G_1(x,s)G_1(y,x)G_2(y,s) U(s,y).
\eqno\,(5.1)
$$
We multiply both sides of this relation by $U(x,t)$
and sum it over $x$; then we get
$$
{\bf E}\,U_{G_1}^{\circ }(t,t)G_2(y,y)=vq_2(y){\bf E}\,
U_{G_1}^{\circ}(t)G_2(y,y)\,U_{G_2}^{\circ }(y)+
$$
$$
vq_2(y)\sum_s{\bf E}\,U_{G_1}^{\circ }(t)\left[ G_2(y,s)\right]^2
U(s,y)+
$$
$$
2vq_2(y)\sum_{s}{\bf E}\,G_1(x,s)G_1(y,x)G_2(y,s) U(s,y) U(x,t).
\eqno\,(5.2)
$$

Regarding $G_1(y,\cdot ) U(\cdot ,t)$ and $G_2(y,\cdot ) U(\cdot ,y)$
in the last term as vectors in $N$-dimensional space, we derive from
estimate (3.19) that
$$
\left| \sum_{s,x}{\bf E}\,G_1(x,s)G_1(y,x)G_2(y,s) U(s,y) U(x,t)
\right| \le
$$
$$
\Vert G_1 \Vert \left( \sum_x \left| G_1(y,x) U(x,t )\right|^2 \right)
^{1/2}\left( \sum_s\left| G_2(y,s) U(s,y)\right|^2 \right) ^{1/2}.
\eqno (5.3)
$$
Inequality (3.18) implies that
the right-hand side of (5.3) is bounded by $b^{-2} \eta^{-3}$.

Let us multiply both sides of (5.2) by $U(y,r)$ and sum them over $y$.
Then one obtains  a relation that together with (3.18) and (5.3)
implies the following estimate for
variable $M_{12}=\sup _x\left( {\bf E}\left| \,U_{G_1}^{\circ }(x)\right|
^2\right) ^{1/2}$:
$$
M_{12}\le v\eta ^{-2}M_{12}+v\eta ^{-3}b^{-1}\sqrt{M_{12}}+2v\eta
^{-4}b^{-2}.
$$
This proves (4.18).

Now (4.17)  follows from (4.18) and relation (5.1).

To derive estimate (4.19), let us consider
the variable
$$
{\bf E}\,U_{G_1}^{\circ }(x_1)U_{G_2}^{\circ }(x_2)U_{G_3}^{\circ
}(x_3)U_{G_4}^{\circ }(x_4)=
{\bf E}\,\left[ U_{G_1}^{\circ }(x)U_{G_2}^{\circ
}(x)U_{G_3}^{\circ }(x_3)\right] ^{\circ } U_{G_4}(x_4).
$$
Let is denote $ T = U_{G_1}^{\circ }U_{G_2}^{\circ}U_{G_3}^{\circ }$ and
and $M (x_1,x_2,x_3,t)= \E T^\circ G_4(t,t)$.
We
apply to $G_4(t,t)$ resolvent identity (3.14)
and  obtain relation
$$
\E T^\circ G_4(t,t) = v\zeta_4 \E T^\circ G_4(t,t) U_{G_4}(t)+
$$
$$
v\zeta_4 \E T^\circ \sum_s [G_4(s,t)]^2 U(s,t) +
$$
$$
v\zeta_4 \sum_{(i,j,k)}\E U^\circ_{G_i}(x_i) U^\circ_{G_j}(x_j)
\sum_{x,s,t} G_k(y,s) G_k(t,y) U(y,x_k) G_4(t,s) U(s,t).
\eqno (5.4)
$$
Repeating previous computations and applying similar estimates,
we obtain inequality
$$
\vert \sum _t M(x_1,x_2,x_3,t) U(t,x_4) \vert \le
{v\over \eta} \E \vert T U_{G_4}^\circ(x_4)\vert +
{v\over \eta} \E \vert T \vert \, \E \vert U_{G_4}^\circ(x_4)\vert +
$$
$$
{v\over \eta^3 b} \E \vert T \vert +
{3v\over \eta b^2} \E \vert U_{G_i}^\circ(x_i) U_{G_j}^\circ(x_j)\vert.
\eqno (5.5)
$$
Here we have applied inequalities  (3.18) and (5.3) to
the last two terms of relation (5.4).
Now it is clear that (5.5) implies (4.19).
Proposition 4.3
is proved.$\Box $

\vskip 0.5cm

{\it Proof of Proposition 4.4}.

Estimate (4.20)  follows from relation (5.2) and estimate (4.18).
Regarding (5.1) and summing it  over $x$,
one can easily derive (4.21) with the help of the
arguments used to prove (4.18).

Let us turn to the proof of (4.22). To do this, let us consider the variable
$$
K(x,y)={\bf E}\left\langle G^{\circ }\right\rangle G^{\circ }(x,x)U_G^{\circ
}(y)={\bf E}\left[ \left\langle G^{\circ }\right\rangle U_G^{\circ
}(y)\right] ^{\circ }G(x,x)
$$
and apply to the last expression resolvent identity (3.12) and formulas
(3.13) and (3.16).
We obtain equality that can be written in the following
form with denotation $R=\left\langle G^{\circ }\right\rangle U_G^{\circ }(y)$
$$
{\bf E}\,R^{\circ }G(x,x)=\zeta v{\bf E}\,R^{\circ
}G(x,x)U_G(x)+\sum_{i=1,2,3}\kappa_i(x,y),
\eqno\,(5.6)
$$
where
$$
\kappa_1(x,y)=\zeta v\sum_s{\bf E}\,R^{\circ }G(x,s)G(x,s) U(s,x),
$$
$$
\kappa_2(x,y)=2\zeta v\sum_{s,t}{\bf E}\left\langle G^{\circ }\right\rangle
G(t,s)G(x,t)u_b^2(t,y)G(x,s) U(s,x),
$$
and
$$
\kappa_3(x,y)=2\zeta vN^{-1}\sum_{s,t}{\bf E}G(t,s)G(x,t)U_G^{\circ
}(y)G(x,s)u_b^2(s,x).
$$
Let us use identity
$$
{\bf E}\,R^{\circ }XY={\bf E}\,R X^\circ \,{\bf E}\,Y+{\bf E}\,RY^\circ
{\bf E}\,X+{\bf E}\,RX^{\circ }Y^{\circ }-{\bf E}\,R\,{\bf E}\,X^{\circ
}Y^{\circ }.
$$
and can rewrite (5.6) in the form
$$
{\bf E}\,R^{\circ }G(x,x)=\frac{vq(x)}{1-vq(x)g(x) }\left[ {\bf E}
\,R U_G^{\circ }(x)G^{\circ }(x,x)-{\bf E}\left\langle G^{\circ
}\right\rangle U_G^{\circ }(y){\bf E}G^{\circ }(x)U_G^{\circ }(x)\right] +
$$
$$
\frac{vq(x)}{1-vq(x)g(x)}\sum_{i=1,2,3}\kappa_i(x,y).
\eqno
\,(5.7)
$$

Taking into account relation (4.18), inequalities (3.18) and
(5.3),  we obtain that
$$
\left| \kappa_i(x,y)\right| \le 2\eta ^{-2}b^{-2}\left( {\bf Var}
\,\left\langle G\right\rangle ^{}\right) ^{1/2}\;{\hbox{for}}\;i=1,2
$$
and
$$
\left| \kappa_3(x,y)\right| \le 2\eta ^{-3}b^{-2}N^{-1}.
$$
Using them, we derive from (5.7) inequality
$$
\left| K(x,y)\right| \le 2\eta ^{-1}\left( {\bf Var}\,\left\langle
G\right\rangle \right) ^{1/2}\left\{ \left( {\bf E}\left| U_G^{\circ
}(x)\right| ^4\right) ^{1/2}+b^{-2}\left( {\bf E}\left| U_G^{\circ
}(x)\right| ^2\right) ^{1/2}\right\} +
$$
$$
2\eta ^{-1}b^{-2}\left( {\bf Var}\,\left\langle G\right\rangle ^{}\right)
^{1/2}+2\eta ^{-2}b^{-2}N^{-1}.
$$
This leads to estimate (4.22). Proposition 4.4 is proved.$\Box $

\vskip 0.5cm

{\it Proof of Proposition 4.5.}

This proof of the estimate (4.23) is the most cumbersome.
Here we  have to use
the resolvent identity (3.12) twice. However,
the computations are based on the same inequalities
as those of the proofs of Propositions 4.3 and 4.4.
Therefore we just indicate the principal lines of the proof
and do not present the  derivations of estimates.

To compute the mathematical expectation
$$
{\bf E}\,M(x,s)={\bf E} \left\langle G_1^{\circ }\right\rangle
\left[ G_2(x,s)\right] ^2 ,
$$
let us apply to $G_2(x,s)$ the resolvent identity (3.12). We obtain
equality
$$
{\bf E}\,M(x,s)=\zeta _2{u (0)\over b}
{\bf E}\left\langle G_1^{\circ
}\right\rangle G_2(x,x)-
$$
$$
\zeta _2{\bf E}\left\langle G_1^{\circ }\right\rangle
\sum_{t}G_2(x,s) G_2(x,t)a(t,s) \sqrt{U(t,s)} .
\eqno\,(5.8)
$$
Relation  (4.21) implies that  the first term of the right-hand side of (5.8)
is the value of the order indicated in (4.23).
Let us consider the second term of (5.8). We
compute mathematical expectation with the help of relations (3.13) and
(3.16) and obtain expression
$$
\zeta_2 {\bf E}\left\langle G_1^{\circ }\right\rangle
\sum_{t}G_2(x,s)  G_2(x,t)a(t,s) \sqrt {U(t,s)}=
\sum_{i=1}^5 \Theta_i(x,s),
\eqno (5.9)
$$
where
$$
\Theta _1(x,s)=v\zeta_2 {\bf E}\left\langle G_1^{\circ }\right\rangle
G_2(x,s)G_2(x,s) \,  \E U_G(s) ,
$$
$$
\Theta_2 (x,s)= v\zeta_2 {\bf E}\left\langle G_1^{\circ }\right\rangle
 G_2(x,s)G_2(x,s)   U_G^\circ(s),
$$
$$
\Theta _3(x,s) ={2v\zeta_2 \over N}{\bf E}
\sum_{t} G_1^2(s,t) U (t,s)G_2(x,s)G_2(x,t),
$$
$$
\Theta _4 (x,s)=v\zeta_2  {\bf E}
\left\langle G_1^{\circ }\right\rangle
\sum_{t} \left[
G_2(x,t)\right] ^2 U(t,s) G_2(s,s) ,
$$
and
$$
\Theta_5(x,s)=2v\zeta_2  {\bf E}\,\left\langle G_1^{\circ }\right\rangle
\sum_{t} G_2(x,s)  G_2(x,t)G_2(s,t) U(t,s) .
$$
$\Theta_1$ is of the form $v \zeta_2  \E M(x,s) \E U_G(s)  $ and
can be put to the right-hand side of (5.9).
The terms $\Theta_2$ and $\Theta_3$ are of the order indicated in
the right-hand side of (4.23). This can be shown with the help
of estimates of the form (5.3).

Regarding $\Theta_4$, we apply the resolvent identity
(3.12) to factor  $G_2(s,s)$.
Repeating the usual computations based on (3.13) and (3.16), we obtain
that
$$
\Theta_4(x,s) = v \zeta_2^2 \sum_{t} \E M(x,t) U(t,s) +
v\zeta_2 \Theta_4(x,s) \E U_{G_2}(s) + \Omega(x,s),
\eqno (5.10)
$$
where $\Omega$  gathers the terms that are all of the order
indicated in (4.23).  This can be verified by direct computation
with the use of estimates (4.18), (4.21), and (4.22).
Not to overload this paper, we do not write down the terms
constituting $\Omega$ and do not present their estimates as well.
Relation (5.10) is of the from that leads to the estimates needed
for $\sum_s \E M(x,s) U(s,x)$.

Regarding $\Theta_5(x,s)$,
we apply (3.12) to $G_2(s,t)$ and obtain, after the use of (3.13) and (3.16)
that
$$
\Theta_5(x,s) = 2v\zeta_2^2{u(0)\over b} \E M(x,s) +
v\zeta\Theta_5(x,s)\E U_{G_2}(s) + \Omega'(x,s),
\eqno (5.11)
$$
where $\Omega'(x,s)$ consists of the terms that are
of the order indicated in (4.23). The form of (5.11)
is also such that, being substituted into (5.9) and then into (5.8), it
leads to the estimates needed.
This observations show that (4.23) is true.

\vskip 0.5cm

{\it Proof of Proposition 4.7.}

We prove relation (4.27) with $k=1$ because the general case does not differ
from this one.
To derive relations for the average value of variable
$t_{12}(x,y)={\bf E}G_1(x,y)G_2(x,y)$, we use
identities (3.12)-(3.14) and repeat  the
 proof of Proposition 4.6.
 Simple
computations lead us to equality
$$
t_{12}(x,y)=g_1(x)\zeta _2\delta _{xy}+\zeta
_2v^2t_{12}(x,y)U_{g_2}(y)+
$$
$$
\zeta _2v^2\sum_st_{12}(x,s) U(s,y)\,g_1(y)+
\zeta _2v^2
\sum_{j=1}^4 \Upsilon_j(x,y),
\eqno\,(5.12)
$$
where
$$
\Upsilon_1(x,y)={\bf E}\sum_sG_1(x,y)G_2(x,s)G_2(y,s) U(s,x),
$$
$$
\Upsilon_2(x,y)={\bf E}\sum_sG_1(x,y)G_1(s,y)G_2(x,s) U(s,x),
$$
$$
\Upsilon_3(x,y)={\bf E}G_1(x,y)G_2(x,y)U_{G_2}^{\circ }(y),
$$
and
$$
\Upsilon_4(x,y)={\bf E}G_2^{\circ }(y,y)\sum_sG_1(x,s)G_2(x,s) U(s,y).
$$

It is easy to see that inequality (4.16) implies estimates
$$
\sup_{x,y}\left| \Upsilon_1(x,y)\right| \le b^{-1}\eta ^{-3},
\quad
\sup_x \vert \sum_y \Upsilon_1(x,y) \vert \le b^{-1}\eta ^{-3}.
$$
The same is valid for $\Upsilon_2(x,y)$.
Similar estimates for $\Upsilon_3(x,y)$ and $\Upsilon_4(x,y) $
follow from relations (4.17) and (4.18).
Thus, (5.12) implies that
$$
t_{12}(x,y)=g_1(x)q_2(x)\delta _{xy}+vg_1(y)q_2(y)\left[ t_{12}U\right]
(x,y)+\Delta (x,y),
\eqno\,(5.13)
$$
where
$$
\sup_{x,y}\left| \Delta (x,y)\right| =o(1) \quad {\hbox {and}}\quad
\sup_x \left| \sum_y \Delta (x,y)\right| =o(1)
\eqno (5.14)
$$
in the limit $n,b\to\infty$ (2.12).

We rewrite  relation (5.13) in the matrix form (cf. (4.35))
$$
t_{12}=\left( I-W^{(g,q)}\right) ^{-1}\left[ {\hbox{Diag}}(g_1q_2)+\Delta \right]
= \sum_{m=0}^\infty [W^{(g,q)}]^m \left( {\hbox{Diag}}(g_1q_2) + \Delta\right).
\eqno (5.15)
$$
Now we can apply to (5.15) the same arguments as to (4.35).
Replacing $g_1(x)$ and $q_2(x)$ by $w_1$ and $w_2$, respectively,
we derive from
(5.14) that for $x\in B_{L+Q}$
$$
t_{12}(x,s) =\sum_{m=0}^M\left( w_1w_2\right) ^{m+1}\left[
U^{m}\right] (x,s)+o(1), \,\, n,b\to\infty.
\eqno\,(5.16)
$$
Multiplying both sides  of (5.16) by $U(s,x)$
and summing it over $s$, we
obtain relation
$$
\sum_{\vert s\vert \le n} t_{12}(x,s) U(s,x) =
\sum_{m=0}^M\left( w_1w_2\right) ^{m+1}\left[
U^{m+1}\right] (x,x)+o(1), \,\, n,b\to\infty.
\eqno (5.17)
$$
Now  convergence (4.40) implies
relation that leads, with $M$ replaced by $\infty$, to (4.27).

\vskip 0.5cm

To prove (4.28), let us sum (5.16) over $s$.
The second part of (5.14) tells us that the terms $\Delta$
remains small when summed over $s$.
Thus we can write relations
$$
\sum_s t_{12}(x,s) =\sum_{m=0}^M\left( w_1w_2\right) ^{m+1}\sum_s \left[
U^{m}\right] (x,s)+o(1), \,\, n,b\to\infty.
\eqno\,(5.18)
$$
Taking into account estimates for terms (4.11) and (4.12),
it is easy to observe that convergence
(4.39) together with (5.18) implies (4.28).
$\Box $

\vskip 0.5cm

%section 6%
%file bs6.tex%

%%%%%%%%%%%%%%%%%%%%%%%%%%%%%%%%%%%%%%%%%%%%%%%%%%%%%%%%%%%%%%%%%%%%%%%

%%%%%%%%%%%%%%%%%%%%%%%%%%%%%%%%%%%%%%%%%%%%%%%%%%%%%%%%%%%%%%%%%%%%%%%

%%%%%%%%%%%%%%%%%%%%%%%%%%%%%%%%%%%%%%%%%%%%%%%%%%%%%%%%%%%%%%%%%%%%%%%

%%%%%%%%%%%%%%%%%%%%%%%%%%%%%%%%%%%%%%%%%%%%%%%%%%%%%%%%%%%%%%%%%%%%%%%

%%%%%%%%%%%%%%%%%%%%%%%%%%%%%%%%%%%%%%%%%%%%%%%%%%%%%%%%%%%%%%%%%%%%%%%

%%%%%%%%%%%%%%%%%%%%%%%%%%%%%%%%%%%%%%%%%%%%%%%%%%%%%%%%%%%%%%%%%%%%%%%

%%%%%%%%%%%%%%%%%%%%%%%%%%%%%%%%%%%%%%%%%%%%%%%%%%%%%%%%%%%%%%%%%%%%%%%

%%%%%%%%%%%%%%%%%%%%%%%%%%%%%%%%%%%%%%%%%%%%%%%%%%%%%%%%%%%%%%%%%%%%%%%

%%%%%%%%%%%%%%%%%%%%%%%%%%%%%%%%%%%%%%%%%%%%%%%%%%%%%%%%%%%%%%%%%%%%%%%

%%%%%%%%%%%%%%%%%%%%%%%%%%%%%%%%%%%%%%%%%%%%%%%%%%%%%%%%%%%%%%%%%%%%%%%

\section{Asymptotic properties of $S(z_1,z_2)$}

In the last decade, the main focus of
the spectral theory of random matrices is
related with the universality conjecture
of local spectral statistics
put forward first by F. Dyson \cite{D}.
This problem is addressed in a large number
of papers where various random matrix ensembles
are studied using different approaches
(see e.g. the review \cite{MGZ}).
The best understood are the Gaussian Unitary Ensemble (GUE)
and its real symmetric analogue GOE (see (2.5)).
The probability distribution of these ensembles
are invariant with respect to the unitary
(orthogonal) transformations.
This leads to the fact that
the joint probability distribution
of eigenvalues of these ensembles
does not depend on the distribution of eigenvectors
and
is given in explicit form \cite{M}.
This allows one to use the powerful technique of the
orthogonal polynomials that provides a detailed information
of the spectral properties of GUE and GOE and related
ensembles on the local scale (see \cite{B,TW} for the
initial results for Gaussian ensembles and \cite{BI,DIZ}
for their generalizations).

The case of band random matrices is different
because the probability distribution of the ensemble
$H^{(n,b)}$ (2.4)
is no more invariant under transformations of the coordinates.
One of the possible ways to study the spectral properties of $H^{(n,b)}$
is to follow
the resolvent expansions
approach
well-known in theoretical physics (see, for example \cite{FMP}).
A rigorous version of it has been
developed
in a series of papers \cite{KP1,KP2,KKP}.

In frameworks of the resolvent approach
(see \cite{KKP} for details),
one considers
 the correlation function
$
C_{n,b}(z_1,z_2), \I z_j \neq 0
$ (3.4)
in the limit when the dimension of the matrix $N$ infinitely increases.
Asymptotic expression for $S(z_1,z_2)$
regarded in the limit
$
z_1 = \lambda_1 + {\hbox{i}}0,\,\, z_2 = \lambda_2 - {\hbox{i}}0
$
supplies one with the information
about  the local properties of eigenvalue distribution  
provided 
$\lambda_1-\lambda_2= O(N^{-1})$. Indeed, according to (2.7), the formal definition
of the eigenvalue density $\rho_{n,b}(\l) = \sigma_{n,b}'(\l)$ is
$$
\rho_{n,b}(\lambda) = {1\over 2\i}\left[ f_{n,b}(\l +\i 0) - f_{n,b}(\l-\i 0)\right].
$$
Then one can consider
expression
$$
R_{n,b}(\l_1,\l_2)=-{1\over 4} \sum_{\delta_1,\delta_2 = -1,+1}
\delta_1 \delta_2 C_{n,b}(\l_1+ \i \delta_1 0, \l_2+ \i \delta_2  0)
$$
as the correlation function of $\rho_{n,b}$.
In general, even if $R_{n,b}$ can be rigorously determined,
it is difficult to carry out the direct study of it.
Taking into account relation (3.5),
one can pass to
more simple expression
$$
\Sigma_{n,b}(\l_1,\l_2) = -{1\over 4 Nb }\sum_{\delta_1,\delta_2 = -1,+1}
\delta_1 \delta_2 S(\l_1+ \i \delta_1 0, \l_2+\i \delta_2 0)
\eqno (6.1)
$$
 and assume that it
corresponds to the leading term of $R_{n,b}(\l_1,\l_2)$  in the limit $n,b\to\infty$.

In present section we follow the same heuristic scheme.
It should be noted that for Wigner random matrices
this approach is justified by
the study of the simultaneous limiting
transition $N\to\infty, \I z_j \to 0$ in the
studies of  $C_{N}(z_1,z_2)$
\cite{BK,K2}.

\vskip 0.5cm
\newpage

{\bf Theorem 6.1}.

 {\it Let $S(z_1,z_2)$ is given by (3.6). Assume that function }$
\tilde u_F(p)${\it \ is such that there exist positive constants
$c_1, \delta$
and $\nu>1$
that
$$
\tilde u_F(p)=\tilde u_F(0)-c_1\left| p\right| ^\nu + o(\left| p\right| ^\nu)
\eqno\,(6.2)
$$
for all $p$ such that  $\vert p\vert \le \delta, \, \delta \to 0$.
Then }
$$
\Sigma_{n,b}(\l_1,\l_2) =
 {1\over Nb} \frac{c_2}{\vert\l_1-
\l_2\vert^{2-1/\nu }}
\left( 1+o(1)\right)
\eqno \,(6.3)
$$
{\it for  $\l_j, j=1,2$ satisfying}
$$
\l_1,\l_2 \to \l \in (-2\sqrt v, 2\sqrt v).
\eqno (6.4)
$$

\vskip 0.5cm

{\it Proof of Theorem 6.1.}

Let us start with the terms of (6.1) that correspond to $\delta_1\delta_2 = -1$.
It follows from (2.9) that
$$
1-vw_1w_2=\frac{z_1-z_2}{w_1-w_2}.
\eqno (6.5)
$$
Also for the real and imaginary parts of
$w(\l+\i 0) = \tau(\l) + \i \rho(\l)$, we have
$$
\tau^2 = {\l^2\over 4v^2},\quad \rho^2 = {4v - \l^2\over 4v^2}
\eqno (6.6)
$$
(here and below we  omit
variable $\l$).
This implies existence of the limits $w(z_1) = \overline {w(z_2)}$  for (6.4).
One can easily deduce from (6.5) that in the limit (6.4) 
$$
1 - vw(z_1) w(z_2) = {\l_1 - \l_2\over 2\i \rho(\l)} = o(1).
\eqno (6.7)
$$
Also we have that
$$
(1-vw_1^2)(1-vw_2^2) = 2 - 2v(\tau^2-\rho^2) = 4v \rho^2.
\eqno (6.8)
$$

Now let us consider $Q(z_1,z_2)$ (3.8) and write that
$$
Q(z_1,z_2) =
{1\over 2\pi} \left(
\int_{-\delta }^\delta +\int_{{\bf R}
\,\setminus (-\delta ,\delta )}
\right)
\frac{w_1^2 w_2^2 \tilde u_F(p)}
{\left[ 1-vw_1w_2\tilde u_F(p)\right] ^2}\d p = {\bf I}_1 + {\bf I_2}.
$$
Relations (6.5) and (6.7) imply equality (cf. (3.9))
$$
\left[ 1-vw_1w_2\tilde u_F(p)\right] ^2 = [\tilde u_F(p) - 1]^2(1+o(1)).
\eqno (6.9)
$$
Since $u(t)$ is monotone, then
$$
{\hbox{liminf}}_{p\in {\bf R}\setminus (-\delta,\delta)} [\tilde u_F(p) - 1]^2 >0.
$$
This means that ${\bf I}_2 < \infty$ in the limit (6.4).

Regarding (6.7), we can write that
in the limit (6.4)
$$
{\bf I}_1 =
\int_{-\delta }^\delta \frac{(2\pi)^{-1}w_1^2 w_2^2 \tilde u_F(p) }
{\left( 1-vw_1w_2^{\ }+ vw_1w_2 \left[\tilde u_F(p) - 1\right] \right) ^2}\d p =
\int_{-\delta }^\delta \frac{(2\pi v)^{-1} \tilde u_F(p) (1+o(1))}
{\left( {z_1-z_2\over w_1 - w_2} + [\tilde u_F(p) - 1] \right)^2}\d p.
$$
Then we derive relation
$$
{\bf I}_1(\l_1+\i 0,\l_2 - \i 0 ) +
{\bf I}_1(\l_1-\i 0,\l_2 + \i 0 ) =
$$
$$
{1\over \pi}\int_{-\delta}^\delta
{
[\tilde u_F(p) - 1]^2 - \left( {\l_1-\l_2\over 2\rho}\right)^2
\over
\left[ [\tilde u_F(p) - 1]^2 + \left( {\l_1-\l_2\over 2\rho}\right)^2\right]^2
}\tilde u_F(p) (1+o(1))
\d p,
\eqno (6.10)
$$
where $o(1)$ corresponds to (6.9) regarded in the limit (6.4).

Now let us use condition (6.2) and
observe that
$$
{1\over \pi} \int_{-\delta} ^{\delta}
{
c_1^2 p^{2\nu} + o(p^{2\nu}) - D^2
\over
\left[c_1^2 p^{2\nu} + o(p^{2\nu}) + D^2\right]^2
}
\d p =
{2\over \pi D^{2-1/\nu}} \int_0^{\delta D^{-1/\nu}}
{
c_1^2 s^{2\nu} + o(s^{2\nu}) -1
\over \left[c_1^2 s^{2\nu} + o(s^{2\nu}) + 1\right]^2
}
\d s,
$$
where we denoted $D = \vert \l_1 - \l_2\vert /(2\rho)$ and $o(p^{2\nu})$
corresponds to the limit $\delta\to 0$ (6.2).
Now it is clear that if we take $\delta$ such that
$\delta \vert \l_1 - \l_2\vert^{-1/\nu} \to \infty$, we
obtain asymptotically
$$
I_1 + \bar I_1 = 4B_\nu(c_1){ (2\rho)^{2-1/\nu} \over \vert \l_1-\l_2\vert^{2-1/\nu} }
\eqno (6.11)
$$
where
$$
B_\nu (c_1) ={1\over 2\pi c_1^{1/\nu}}\left[
\int_{0}^\infty {\d s\over 1+ s^{2\nu}  }
-
2\int_{0}^\infty {\d s\over (1+ s^{2\nu})^2  }\right].
\eqno (6.12)
$$

To prove relation (6.3), it remains to consider the sum
$$
I(\l_1+\i 0, \l_2+\i 0 ) +
I(\l_1-\i 0, \l_2-\i 0 ).
$$
It is easy to observe that
 relations of the form (6.8) imply boundedness of this sum in the limit
 (6.4)

Now gathering relations (6.8) and (6.11), we derive that
$$
\Sigma_{n,b} (\l_1,\l_2) s=
{1\over Nb}
{B_\nu(c_1) \over (2\rho)^{1/\nu}}
{1 \over
\vert \l_1-\l_2\vert^{2-1/\nu}
} (1+o(1))
\eqno (6.13)
$$
This proves  (6.3).$\Box$

\vskip 0.5cm

Let us discuss two  consequences of Theorem 6.1.
Let us assume first that $u(t)$ is such that
$$
u_2\equiv \int t^2  u(t)\,\d t<\infty .
\eqno
\,(6.14)
$$
Then (6.2) holds with $\nu =2$ and $c_1 = u_2$.
Regarding the right-hand side of (3.5)
in the limit (6.4) with $\l_j =\l + r_j/N$, $j=1,2$,
we obtain asymptotic relation
$$
\Sigma(\l_1,\l_2)={\sqrt{N} \over b}
{B_2(u_2) \over 2 (2\rho)^{1/2}}
{1 \over \vert r_1-r_2\vert ^{3/2}}
(1+o(1)),
\eqno
\,(6.15a)
$$
where
$$
B_2 (u_2)= -  {1\over 4\pi  \sqrt{u_2} } \int_0^\infty {\d s\over 1+s^4} =
 - {1\over 4\pi \sqrt{u_2} } \Gamma\left(5\over 4\right)  \Gamma\left(3\over
4\right).
\eqno (6.15b)
$$

\vskip 0.5cm
Now let us assume that (6.14) is not true. Suppose that there exists such
$
1<\nu' <2$ that
$$
 u(t)=O(\left| t\right| ^{-1-\nu'})
\quad {\hbox{as \, }}\;t\rightarrow \infty .
\eqno
\,(6.16)
$$
Then one can easily derive that (6.2) holds with $\nu=\nu'$. This follows from
elementary computations based on equalities
$$
\tilde u_F(p)= \tilde u_F(0) -\int_{-\infty }^\infty \left( 1-\cos pt\right)
u(t)\d t
$$
and
$$
\frac 1p\int_{-\infty }^\infty \left( 1-\cos y\right) u(yp^{-1})\d y =
\left| p\right| ^{\nu'} \int_{-\infty }^\infty \frac{\left(
1-\cos y\right) }{\left| y\right| ^{1+\nu' }}\d y + o(\vert p\vert^{\nu'}),
\,\,p\to 0.
$$
Therefore, if (6.16) holds, then
$$
\Sigma(\l_1,\l_2)= 
{N^{1-1/\nu } \over b}
{B_\nu(c_1) \over (2\rho)^{1/\nu}}
{1 \over \vert r_1-r_2\vert^{2-1/\nu }  }\left( 1+o(1)\right) .
\eqno
\,(6.17)
$$

The form of 
asymptotic expressions (6.15a) and (6.17)
coincides with that determined by Altshuler and Shklovski
for  the spectral correlation function
of band random matrices  (see \cite{MFD}
for this and similar results).
In these works, the factor
$\vert r_1 - r_2 \vert ^{-3/2}$ appeared
instead of usual for random matrices expression
$\vert r_1-r_2\vert^{-2}$ (see (3.10)). This
has been interpreted as the evidence  of (relatively)
localized eigenvectors of $H^{(n,b)}$ in the limit
$1\ll b\ll n $ with the localization
length $b^2/n$.
Let us note that the asymptotic expressions similar to  (6.15) have also appeared
in the recent work \cite{Syl}, where the band random matrix ensemble
$H^{(n,b)}$ was considered under condition (6.14).
However it should be stressed that no explicit
expressions like  (3.6)  and (6.15) were obtained neither in
\cite{MFD} nor in \cite{Syl}.

\section{Summary}

We consider a family of random matrix ensembles $\left\{ H^{(n,b)}\right\} $
of the band-type form. More precisely, we are related with real symmetric
$N\times N$ matrices, $N=2n+1$, whose entries are jointly independent Gaussian random
variables with zero mean value.
The band-type form means that the variance
of the matrix entries $H^{(n,b)}(x,y)$ is proportional to
$ u\left(\frac{x-y}b\right)\ge 0 $.

We study asymptotic behavior of the 
correlation function
$$
C_{n,b}(z_1,z_2) = \E f_{n,b}(z_1)  f_{n,b}(z_2) - \E f_{n,b}(z_1)\E f_{n,b}(z_1),
$$
where $ f_{n,b}(z)$ is the  
normalized trace $\left\langle
G^{(n,b)}(z)\right\rangle $ of the resolvent of $H^{(n,b)}$.

We have proved that if $\I z_j$ is large enough, then
in the limit $1\ll b\ll N^{1/3}$ 
$$
C_{n,b}(z_1,z_2)=\frac 1{Nb}S(z_1,z_2)+o\left( \frac 1{Nb}\right).
$$
We have found explicit form of the leading term $S(z_1,z_2)$ in this limit.
Assuming that expression
$
\Sigma_{n,b}(\l_1,\l_2) $ (6.1)
is closely related with the correlation function of the eigenvalue density,
we have studied it in the limit $N,b\to\infty$ and $\l_1 - \l_2 = (r_1-r_2)/N$.

 Our main conclusion is that the limiting expression
for $\Sigma_{n,b}$ exhibits different behavior depending on
the rate of decay of $u(t)$ at infinity.

If $\int t^2 u(t) \d t <\infty$, then (6.1) is given by 
$$
- C \frac{\sqrt{N}}b\frac 1{\vert r_1-r_2\vert^{3/2}}(1+o(1)), \,\, C>0.
$$
If $\tilde u(t)=O(\left| t\right| ^{-1-\nu })$ with $1<\nu <2$, then the
asymptotic expression for (7.1) is proportional to 
$$
\frac{n^{1-1/\nu }}b\frac 1{\vert r_1-r_2\vert^{2-1/\nu }}.
$$
In both cases the exponents do not depend on the particular form
of the function $u(t)$. Moreover, in the first case
the exponents do not depend on $u$ at all. This can be regarded
as a kind of spectral universality for band random matrix ensembles.
On can conject that these characteristics also do not depend
on the probability distribution of the random variables $a(x,y)$ (2.1).

Our results show that $S(z_1,z_2)$ determines at least two scales of
universality
in the local spectral properties of band-type random matrices. These scales
coincide with those detected in theoretical physics  for the
(relative) localization length and density-density correlation function for
these ensembles \cite{MFD}. In the papers cited also the third scale when 
$ u(t)=O(\left|
t\right| ^{-\gamma })$ with $\gamma \in (1/2,1)$ has been observed.
It  been shown to produce the asymptotics $N^{-2}\vert r_1-r_2\vert^{-2}$
which is typical for ''full'' random matrices like GOE
\cite{B,D,FMP,KKP}. Unfortunately, this
asymptotic regime for band random matrices
is out of reach of our technique.

\vskip 0.5cm

{\bf Acknowledgments}.

The first author is grateful to Ya. Fyodorov for fruitful discussions and
explanations of the role of the Altshuler-Shklovski asymptotics.
The financial support from SFB237 (Germany) during autumns 1997 and 1999 and
kind hospitality of Ruhr-University Bochum is gratefully acknowledged by the
first author.
We also thank the anonymous
referee for useful remarks and conjectures concerning
the technical questions and general exposition as well.

\end{document}